
 \def\LOADED{\relax}

 \ifx\BoxedEPSFLoaded\LOADED
  \immediate\write16{}
  \immediate\write16{  BoxedEPSF macros already defined!}
  \immediate\write16{}
   
 \fi

 \let\BoxedEPSFLoaded\LOADED

 \chardef\CatAt\the\catcode`\@
 \chardef\CatColon\the\catcode`\:
 \catcode`\@=11
 \catcode`\:=12

 \let\wlog@ld\wlog
 \def\wlog#1{\relax}

 \newif\ifIN@
 \newdimen\XShift@ \newdimen\YShift@
 \newtoks\Realtoks

%
 \newdimen\Wd@ \newdimen\Ht@
 \newdimen\Wd@@ \newdimen\Ht@@
 \newdimen\TT@
 \newdimen\LT@
 \newdimen\BT@
 \newdimen\RT@
 \newdimen\XSlide@ \newdimen\YSlide@
 \newdimen\TheScale 
 \newdimen\FigScale 
 \newdimen\ForcedDim@@

 \newtoks\EPSFDirectorytoks@
 \newtoks\EPSFNametoks@
 \newtoks\BdBoxtoks@

 \newif\ifNotIn@
 \newif\ifForcedDim@
 \newif\ifForcedHeight@
 \newif\ifPSOrigin

 \newread\EPSFile@

 \newif\ifIN@\def\IN@{\expandafter\INN@\expandafter}
  \long\def\INN@0#1@#2@{\long\def\NI@##1#1##2##3\ENDNI@
    {\ifx\m@rker##2\IN@false\else\IN@true\fi}%
     \expandafter\NI@#2@@#1\m@rker\ENDNI@}
  \def\m@rker{\m@@rker}

  \newtoks\Initialtoks@  \newtoks\Terminaltoks@
  \def\SPLIT@{\expandafter\SPLITT@\expandafter}
  \def\SPLITT@0#1@#2@{\def\TTILPS@##1#1##2@{%
     \Initialtoks@{##1}\Terminaltoks@{##2}}\expandafter\TTILPS@#2@}


  \newtoks\Trimtoks@

 \def\ForeTrim@{\expandafter\ForeTrim@@\expandafter}
 \def\ForePrim@0 #1@{\Trimtoks@{#1}}
 \def\ForeTrim@@0#1@{\IN@0\m@rker. @\m@rker.#1@%
     \ifIN@\ForePrim@0#1@%
     \else\Trimtoks@\expandafter{#1}\fi}

  \def\Trim@0#1@{%
      \ForeTrim@0#1@%
      \IN@0 @\the\Trimtoks@ @%
        \ifIN@
             \SPLIT@0 @\the\Trimtoks@ @\Trimtoks@\Initialtoks@
             \IN@0\the\Terminaltoks@ @ @%
                 \ifIN@
                 \else \Trimtoks@ {FigNameWithSpace}%
                 \fi
        \fi
      }


   \newtoks\pt@ks
   \def \getpt@ks 0.0#1@{\pt@ks{#1}}
   \dimen0=0pt\expandafter\getpt@ks\the\dimen0@

  \newtoks\Realtoks
  \def\Real#1{%
    \dimen2=#1%
      \SPLIT@0\the\pt@ks @\the\dimen2@
       \Realtoks=\Initialtoks@
            }

   \newdimen\Product
   \def\Mult#1#2{%
     \dimen4=#1\relax
     \dimen6=#2%
     \Real{\dimen4}%
     \Product=\the\Realtoks\dimen6%
        }

 \newdimen\Inverse
 \newdimen\hmxdim@ \hmxdim@=8192pt
 \def\Invert#1{%
  \Inverse=\hmxdim@
  \dimen0=#1%
  \divide\Inverse \dimen0%
  \multiply\Inverse 8}

   \def\Rescale#1#2#3{
              \divide #1 by 100\relax
              \dimen2=#3\divide\dimen2 by 100 \Invert{\dimen2}%
              \Mult{#1}{#2}%
              \Mult\Product\Inverse
              #1=\Product}

  \def\Scale#1{\dimen0=\TheScale%
      \divide #1 by  1280
      \divide \dimen0 by 5120%
      \multiply#1 by \dimen0
      \divide#1 by 10  
     }


 \newbox\scrunchbox

 \def\Scrunched#1{{\setbox\scrunchbox\hbox{#1}%
   \wd\scrunchbox=0pt
   \ht\scrunchbox=0pt
   \dp\scrunchbox=0pt
   \box\scrunchbox}}

 \def\Shifted@#1{%
   \vbox {\kern-\YShift@
       \hbox {\kern\XShift@\hbox{#1}\kern-\XShift@}%
           \kern\YShift@}}

  %
 \def\cBoxedEPSF#1{{}\leavevmode
   \ReadNameAndScale@{#1}%
   \ReadEPSFile@ \ReadBdB@x
   \edef\EPSFSpec@{\the\EPSFDirectorytoks@\the\EPSFNametoks@}
     \TrimFigDims@
     \CalculateFigScale@
     \ScaleFigDims@
     \SetInkShift@
     \Real{\FigScale}\edef\FigSc@leReal{\the\Realtoks}%
   \hbox{$\mathsurround=0pt\relax
         \vcenter{\hbox{%
             \FrameSpider{\hskip-.4pt\vrule}%
             \vbox to \Ht@{\parindent=\z@%
                \FrameSpider{\vskip-.4pt\hrule}\vfil
                \hbox to \Wd@{\hfil}%
                \vfil
                \InkShift@{\EPSFSpecial{\EPSFSpec@}{\FigSc@leReal}}%
             \FrameSpider{\hrule\vskip-.4pt}}%
         \FrameSpider{\vrule\hskip-.4pt}}}%
     $}%
    \CleanRegisters@
    }

 \def\tBoxedEPSF#1{\setbox4\hbox{\cBoxedEPSF{#1}}%
     \setbox4\hbox{\raise -\ht4 \hbox{\box4}}%
     \box4
      }

 \def\bBoxedEPSF#1{\setbox4\hbox{\cBoxedEPSF{#1}}%
     \setbox4\hbox{\raise \dp4 \hbox{\box4}}%
     \box4
      }

  \let\BoxedEPSF\cBoxedEPSF

  %

  \def\EPSFxsize{\afterassignment\ForceW@\ForcedDim@@}
      \def\ForceW@{\ForcedDim@true\ForcedHeight@false}

  \def\EPSFysize{\afterassignment\ForceH@\ForcedDim@@}
      \def\ForceH@{\ForcedDim@true\ForcedHeight@true}

%
 \def\ReadNameAndScale@#1{\IN@0 scaled@#1@
   \ifIN@\ReadNameAndScale@@0#1@%
   \else \ReadNameAndScale@@0#1 scaled\DefaultMilScale @
   \fi}

 \def\ReadNameAndScale@@0#1scaled#2@{
    \Trim@0#1@%
    \EPSFNametoks@\expandafter{\the\Trimtoks@}%
    \FigScale=#2 pt%
     }

 \def\SetDefaultEPSFScale#1{%
      \global\def\DefaultMilScale{#1}}

 \SetDefaultEPSFScale{1000}

%

 \def \SetBogusBbox@{%
     \global\BdBoxtoks@{ BoundingBox:0 0 100 100 }%
     \global\def\BdBoxLine@{ BoundingBox:0 0 100 100 }%
     }

 \def\ReadEPSFile@{%
     \openin\EPSFile@
       =\the\EPSFDirectorytoks@\the\EPSFNametoks@
     \relax 
  \ifeof\EPSFile@
     \SetBogusBbox@
     \immediate\write16{}%
     \message{ *** EPS FILE  }%
     \message\expandafter{\the\EPSFNametoks@}%
     \message{ NOT FOUND!  }%
     \immediate\write16{}\relax%
  \else
   \begingroup
   \catcode`\%=12\catcode`\:=12\catcode`\\=12
   \NotIn@true
    \loop
      \ifeof\EPSFile@\NotIn@false
        \SetBogusBbox@
        \immediate\write16{}%
        \message{ *** BoundingBox not found in }%
        \message\expandafter{\the\EPSFNametoks@\space *** }%
        \immediate\write16{}%
      \else\global\read\EPSFile@ to \BdBoxLine@
      \fi
      \global\BdBoxtoks@\expandafter{\BdBoxLine@}%
      \IN@0BoundingBox:@\the\BdBoxtoks@ @%
      \ifIN@\NotIn@false\fi%
    \ifNotIn@\repeat
   \endgroup\relax
  \fi
  \closein\EPSFile@
   }

  \def\ReadBdB@x{
   \expandafter\ReadBdB@x@\the\BdBoxtoks@ @}

  \def\ReadBdB@x@#1BoundingBox:#2@{
    \ForeTrim@0#2@%
    \expandafter\ReadBdB@x@@\the\Trimtoks@ @%
   }

  \newtoks\LLXtoks@ 
  \newtoks\LLYtoks@

  \def\ReadBdB@x@@#1 #2 #3 #4@{
      \Wd@=#3bp\advance\Wd@ by -#1bp%
      \Ht@=#4bp\advance\Ht@ by-#2bp%
       \Wd@@=\Wd@ \Ht@@=\Ht@
       \LLXtoks@={#1}\LLYtoks@={#2}
      \ifPSOrigin\XShift@=-#1bp\YShift@=-#2bp\fi
     }

  %
  \def\SetEPSFDirectory{
           \bgroup\catcode`\:=12\relax
           \def\G@bbl@##1{}\def\\{}
           \global\edef\B@ckslash{\expandafter\G@bbl@\string\\}
           \global\let\bs\B@ckslash\relax
           \SetEPSFDirectory@}

 \def\SetEPSFDirectory@#1{
    \Trim@0#1@
    \global\EPSFDirectorytoks@\expandafter{\the\Trimtoks@ }\relax
    \egroup}

%
 \def\TrimTop#1{\advance\TT@ by #1}
 \def\TrimLeft#1{\advance\LT@ by #1}
 \def\TrimBottom#1{\advance\BT@ by #1}
 \def\TrimRight#1{\advance\RT@ by #1}

 \def\TrimFigDims@{%
    \advance\Wd@ by -\LT@
    \advance\Wd@ by -\RT@ \RT@=\z@
    \advance\Ht@ by -\TT@ \TT@=\z@
    \advance\Ht@ by -\BT@
    }

%
  \def\ForceWidth#1{\ForcedDim@true
       \ForcedDim@@#1\ForcedHeight@false}

  \def\ForceHeight#1{\ForcedDim@true
       \ForcedDim@@=#1\ForcedHeight@true}

  \def\epsfxsize{\afterassignment\ForceW@\ForcedDim@@}
      \def\ForceW@{\ForcedDim@true\ForcedHeight@false}

  \def\epsfysize{\afterassignment\ForceH@\ForcedDim@@}
      \def\ForceH@{\ForcedDim@true\ForcedHeight@true}

  \def\CalculateFigScale@{%
     \ifForcedDim@\FigScale=1000pt
           \ifForcedHeight@
                \Rescale\FigScale\ForcedDim@@\Ht@
           \else
                \Rescale\FigScale\ForcedDim@@\Wd@
           \fi
     \fi}

  \def\ScaleFigDims@{\TheScale=\FigScale
      \ifForcedDim@
           \ifForcedHeight@ \Ht@=\ForcedDim@@  \Scale\Wd@
           \else \Wd@=\ForcedDim@@ \Scale\Ht@
           \fi
      \else \Scale\Wd@\Scale\Ht@
      \fi
      \ForcedDim@false
      \Scale\LT@\Scale\BT@ 
      \Scale\XShift@\Scale\YShift@
      }

 
 \def\ShowReservedBoxes{\gdef\FrameSpider##1{##1}}
 \let\HideDisplacementBoxes\HideReservedBoxes 

  \ShowReservedBoxes

 \def\hSlide#1{\advance\XSlide@ by #1}
 \def\vSlide#1{\advance\YSlide@ by #1}

  \def\SetInkShift@{%
            \advance\XShift@ by -\LT@
            \advance\XShift@ by \XSlide@
            \advance\YShift@ by -\BT@
            \advance\YShift@ by -\YSlide@
             }
  \def\InkShift@#1{\Shifted@{\Scrunched{#1}}}

  %
  \def\CleanRegisters@{%
      \globaldefs=1\relax
        \XShift@=\z@\YShift@=\z@\XSlide@=\z@\YSlide@=\z@
        \TT@=\z@\LT@=\z@\BT@=\z@\RT@=\z@
      \globaldefs=0\relax}


 \def\SetTexturesEPSFSpecial{\PSOriginfalse
  \gdef\EPSFSpecial##1##2{\relax
    \edef\specialthis{##2}%
    \SPLIT@0.@\specialthis.@\relax
    \special{illustration ##1 scaled
                        \the\Initialtoks@}}}

  \def\SetUnixCoopEPSFSpecial{\PSOrigintrue
   \gdef\EPSFSpecial##1##2{%
      \dimen4=##2pt
      \divide\dimen4 by 1000\relax
      \Real{\dimen4}
      \edef\Aux@{\the\Realtoks}%
      \includegraphics{##1}}}

  \def\SetRokickiEPSFSpecial{\PSOrigintrue
   \gdef\EPSFSpecial##1##2{%
      \dimen4=##2pt
      \divide\dimen4 by 10\relax
      \Real{\dimen4}
      \edef\Aux@{\the\Realtoks}%
      \includegraphics{##1}}}

%
\def\SetOzTeXEPSFSpecial{\PSOrigintrue
 \gdef\EPSFSpecial##1##2{%
	\dimen4=##2pt
	\divide\dimen4 by 1000\relax
	\Real{\dimen4}
	\edef\Aux@{\the\Realtoks}
	\special{epsf=##1\space scale=\Aux@}}}


 \def\SetArborEPSFSpecial{\PSOriginfalse
   \gdef\EPSFSpecial##1##2{%
     \edef\specialthis{##2}%
     \SPLIT@0.@\specialthis.@\relax
     \special{ps: epsfile ##1\space \the\Initialtoks@}}}

 \def\SetClarkEPSFSpecial{\PSOriginfalse
   \gdef\EPSFSpecial##1##2{%
     \Rescale {\Wd@@}{##2pt}{1000pt}%
     \Rescale {\Ht@@}{##2pt}{1000pt}%
     \special{dvitops: import
           ##1\space\the\Wd@@\space\the\Ht@@}}}


  \def\SetBeebeEPSFSpecial{
   \PSOriginfalse
   \gdef\EPSFSpecial##1##2{\relax
    \special{language "PostScript",
      literal "##2 1000 div ##2 1000 div scale
       \the\LLXtoks@\space neg \the\LLYtoks@\space neg translate",
              overlay "##1"}}}

 \def\SetStandardEPSFSpecial{%
   \gdef\EPSFSpecial##1##2{%
     \immediate\write16{}
     \immediate\write16{%
       **** Sorry! There is still no standard for \string%
       \special \space EPSF integration *****}%
     \immediate\write16{%
      --- So you will have to identify your driver using a command}%
     \immediate\write16{%
      --- of the form \string\Set...EPSFSpecial, in order to get}%
     \immediate\write16{%
      --- your graphics to print.  See BoxedEPSF.doc.}%
     \immediate\write16{}
     \KillEPSFSpecial
     }}

  \def\KillEPSFSpecial{\gdef\EPSFSpecial##1##2{}}

  \SetStandardEPSFSpecial

 \let\wlog\wlog@ld

 \catcode`\@=\CatAt
 \catcode`\:=\CatColon

\SetRokickiEPSFSpecial  


\HideDisplacementBoxes     

\input amstex
\documentstyle{amsppt}

\newcount\nummereq 
\nummereq=0
\def\labeleq#1{\global\advance\nummereq
by1\global\edef\here{\number\nummereq}\here
\global\let #1 =\here}

\newcount\nummer
\def\newsec#1{\def\curstem{#1} \nummer=0}
\def\numbit{\global\advance\nummer
by1\global\edef\here{\curstem.\number\nummer}\message{
\here}\here}

\magnification =1200
\hsize =15truecm
\hcorrection{.5truein}
\vsize =22.5truecm
\hfuzz1pc               
\TagsOnRight
\NoBlackBoxes
\expandafter\redefine\csname logo\string @\endcsname{}  

\define\C{\Bbb C}
\define\R{\Bbb R}
\define\di{\partial}
\define\restr{\negthickspace \mid}
\redefine\qed{$\square$}
\define\sn{\text{\rm sn\,}}
\define\cn{\text{\rm cn\,}}
\define\dn{\text{\rm dn\,}}
\def\hedd#1{\noindent{\it #1\ }}

\define\BB/{B\"acklund}
\topmatter
\title\nofrills B\"acklund transformations and knots of constant torsion \endtitle
\author Annalisa Calini, Thomas Ivey \endauthor
\address Dept. of Mathematics, Case Western Reserve Univ.,
Cleveland OH 44106\endaddress
\curraddr Dept. of Mathematics, College of Charleston, Charleston SC 29424
\endcurraddr
\email {calini\@math.cofc.edu} \endemail
\address Dept. of Mathematics, Case Western Reserve Univ.,
Cleveland OH 44106\endaddress
\email {txi4\@po.cwru.edu} \endemail
\abstract\nofrills {{\bf Abstract}
The B\"acklund transformation for pseudospherical surfaces, which is 
equivalent to that of the sine-Gordon equation, can be restricted
to give a transformation on space curves that preserves constant torsion.  We
study its effects on closed curves (in particular, elastic rods)
 that generate multiphase solutions for
the vortex filament flow (also known as the Localized  Induction Equation). 
In doing so, we obtain analytic constant-torsion representatives for a
large number of knot types.
}\endabstract
\footnote""{
}
\endtopmatter

\document
\heading Introduction \endheading
Soliton equations have become a familiar presence in the 
differential geometry of curves and surfaces. The description of 
pseudo-spherical surfaces and their asymptotic lines in terms of the
sine-Gordon equation dates back nearly a century. Of roughly
the same date is the derivation (see [Ri]) of the Localized Induction Equation
(LIE), which provides one of the richest examples of connection between
curve geometry and integrability. 

The understanding of this connection has progressed in recent years
along different directions. On the one hand, several fundamental
properties of soliton equations have been given a geometrical 
realization; in the case of the LIE, its bihamiltonian structure and
 recursion operator, its hierarchy of constants of motion, and
its relation to the nonlinear Schr\"{o}dinger equation possess a natural
geometric interpretation [L-P1]. 

On the other hand, some well-known classes of curves in
differential geometry have been identified with solutions of
integrable equations: for example, elastic curves and center-lines of
elastic rods are among the solitons for the LIE, curves of constant
torsion correspond to characteristics for the sine-Gordon
equation, and planar and spherical curves are associated with solutions of
the mKdV hierarchy ([G-P],[L-P2],[D-S]).

A third direction of research ([C],[M-R]) concerns the topological
properties of closed curves that arise as solutions to soliton
equations. A major question is whether the presence of infinitely
many symmetries and the associated sequence of integral invariants
may be related to knot invariants; in fact, one hopes that the knot types 
of the periodic analogues of soliton solutions (the so-called multi-phase 
solutions) can be described using methods of integrable systems, such
as the Floquet spectrum and \BB/ transformations. We remark that
the knot types of closed elastic curves, which turn out 
to generate two-phase solutions of the LIE, 
were classified by Langer and Singer [L-S1] with the aid of 
infinitesimal symmetries related to flows in the LIE hierarchy.

In this work we illustrate these different approaches through the
study of a unifying example. First, we propose a geometric 
realization of the \BB/ transformation for the sine-Gordon equation
in the context of curves of constant torsion, by restricting the
classical \BB/ transformation for pseudo-spherical surfaces [C-T] to
a constant-torsion-preserving transformation between asymptotic lines.
Because it is not necessary to embed the curve in a surface to 
carry out the transformation,
this is also a tool for obtaining new interesting 
curves of constant torsion.

Secondly, we are interested in curves that generate special solutions 
(in particular, solitons and their periodic counterparts) of 
integrable equations and whose analytic expressions and general
properties are well-known. This choice is particularly significant
for closed curves, since multi-phase solutions are dense in the
space of all periodic solutions: it is therefore important to understand
their geometrical and topological properties.
The simplest non-planar constant torsion curves in this class are
the centerlines of elastic rods.  (Elastic rods, which also generate
two-phase solutions of the LIE, are critical curves for the elastic
energy (total squared
curvature), subject to fixed total torsion and fixed length.)  

Finally, we study the topological properties of curves obtained by
means of (single and iterated) \BB/ transformations of constant 
torsion elastic rods.
The behaviour of these transformations is diverse: we show 
that a single \BB/ transformation preserving closure leaves the
class of constant torsion elastic rods invariant, and we compute an interesting
formula relating the linking number of a curve with its transform
and the self-linking of the original curve.
The iterated \BB/ transformation instead allows one to leave the class of elastic
rods and to produce a variety of phenomena such as knotting, self-intersections
and unknotting. Effects of \BB/ transformations on the multi-phase solutions
of the sine-Gordon equation have been discussed in the work by
Ercolani, Forest and McLaughlin [E-F-M], which uses a combination of
Floquet theory for the associated spectral problem and methods of algebraic 
geometry. We present here the first concrete realization of their 
conclusions in a geometric setting, and give simple direct proofs 
for the case of elastic rods of constant torsion.

This article is structured as follows. In order to implement the
\BB/ transformation for pseudospherical surfaces restricted to
constant torsion curves, we need to solve an associated Riccati
equation. To overcome this difficulty, we prove in 
Part I that, if we know the expression of the Frenet frame of a
given curve that depends analytically on the constant torsion
$\tau$, then we can obtain the solution to the Riccati equation
by analytic continuation to an imaginary value of the torsion.
In other words, for any curve whose Frenet 
frame can be continued analytically as a function of $\tau$, with
fixed curvature function, the \BB/ transformation is 
directly computable.  This is illustrated for closed elastic rods
of constant torsion.  Using 
Langer and Singer's solution of the Euler-Lagrange equations [L-S4],
we derive explicit formulas for these curves and their Frenet frames,
and we also describe their knot types.
We then show
that the condition that the \BB/ transformation give another closed curve
forces the new curve to be an elastic rod that is congruent to the
original.
However, such a closed \BB/ transformation carries topological
information by providing a measure of the knottedness of the
initial elastic rod. In \S 1.6 we show this by relating the linking number
of the curve with its transform to the self-linking of the 
original curve.  Part I ends with a discussion of, and 
a precise conjecture about, the effect of \BB/ transformations on 
constant torsion $n$-phase (and $n$-soliton) solutions of the LIE
hierarchy.  The conjecture, which is easily verified for the
first few values of $n$, concerns the form of certain Killing fields
associated with $n$-soliton curves, and it implies that $n$-solitons
are in general taken to $(n+1)$-solitons. 

Part II concerns the iterated \BB/ transformation, which consists of
two successive \BB/ transformations using related solutions of 
the Riccati equations. For the purpose
of constructing iterations, we conveniently rederive the original
\BB/ formula in terms of a gauge transformation of the 
sine-Gordon linear
system. (This is simply the spatial
part, in characteristic coordinates, of the Lax pair,
 and is equivalent to our Riccati equation.)
Using this formulation, a simple algebraic procedure produces 
the iterated formula which we apply to constant torsion elastic 
rods to produce a variety of interesting curves. 
The closure condition, developed using standard arguments of Floquet theory,
 is in this case less restrictive.
As a consequence, we prove that, beyond producing congruent elastic rods, 
an iterated \BB/ transformation can take constant torsion 
two-phase solutions to multi-solitons whose knot type differs
dramatically from that of the original curve. All three
phenomena occur: knotting, self-intersections and unknotting, all
related to the Floquet spectrum of the linear system associated
to the initial curve, as illustrated and discussed in the concluding
section.
\bigskip
The investigations that we report here were initiated by, and advanced with
the steady encouragement of, the geometry group at Case Western Reserve
University (Joel Langer, David Singer, and the authors), with important
contributions from Ron Perline of Drexel University.  In particular,
Perline first asked if there exist constant torsion solitons for the LIE;
Singer classified the elastic rods of constant torsion; and, Langer and Perline
developed the theory of planar-like solitons and their associated Killing fields
which allows us to connect our \BB/ transformations with the LIE hierarchy.
We gratefully acknowledge these contributions, and the lively discussions
that helped to shape this work.

\heading I. The single \BB/ transformation \endheading
\newsec{1}
\heading 1.1 \BB/ transformations and the Frenet frame \endheading
In classical differential geometry, a \BB/ transformation takes a
given pseudospherical (i.e. constant negative Gauss curvature)
surface to a new pseudospherical surface.  As
explained by Chern and Terng [C-T], the new surface is connected to the old surface
by line segments that are tangent to both surfaces, of a fixed length, and such
that the angle between the surface normals at corresponding points is also constant.
Moreover, the \BB/ transformation takes
asymptotic lines to asymptotic lines.  Since the asymptotic lines on a
pseudospherical surface have constant torsion, it is not surprising that we can
restrict the \BB/ transformation to get a transformation that carries constant
torsion curves to constant torsion curves.

\define\gnew{\widetilde{\gamma}}
\define\k{\kappa}
\proclaim{Theorem \numbit} Let $\gamma(s)$ be a smooth curve of constant
torsion $\tau$ in $\R^3$, parametrized by arclength $s$. Let $T,N,B$ be a Frenet frame, and
$\k(s)$ the curvature of $\gamma$.  For any
constant $C$, let $\beta=\beta(s;\k(s),C)$ be a solution of the differential equation
$$\dfrac{d\beta}{ds} = C\sin\beta - \k.\tag\labeleq{\betaODE}$$
Then the curve
$$\gnew(s) = \gamma(s) + \dfrac{2C}{C^2+\tau^2}(\cos\beta T
+\sin\beta N)\tag\labeleq{\curveback}$$
is a curve of constant torsion $\tau$, also parametrized by arclength.
\endproclaim
\demo{Proof}  Since the angle between the normals of the two pseudospherical
surfaces is constant at corresponding points, it follows that the angle
between the binormals of $\gamma$ and $\gnew$ is constant.
In fact, one computes the Frenet frame along $\gnew$:
$$\align
\widetilde T &= T + (1-\cos \theta)\sin\beta(\cos\beta N - \sin \beta T)
+ \sin\theta \sin \beta B \\
\widetilde N &= N - (1-\cos\theta)\cos\beta(\cos\beta N - \sin \beta T)
-\sin\theta \cos\beta B\\
\widetilde B &= \cos\theta B + \sin\theta (\cos\beta N - \sin \beta T),
\endalign$$
where $\tan(\theta/2) = C/\tau$.  Then the torsion is $\tau$ and the curvature
is
\define\knew{\widetilde \k}
$$\knew = \k - 2C \sin \beta.\tag\labeleq{\betak}$$
\qed\enddemo

We should point out that we will assume our curves have a smooth {\it generalized Frenet
frame}, that is, smooth orthonormal vector fields $T,N,B$ that satisfy the Frenet equations
with respect to smooth functions $\k$, $\tau$.
In particular, $\k$ is allowed to change sign.  Given an orientation of
$\gamma$, $T$ is
well defined, and $N,B,\k$ are well-defined
along such a curve up to multiplying all three by $-1$.  If the curve is closed of
length $L$,
then our assumptions imply that there must either be a smooth Frenet
framing by vector fields of period $L$, or a framing by vector fields of period
$2L$ such that $N,B,\k$ change by a minus sign after one circuit of the curve.
(We will refer to these closed curves as `even' and `odd', respectively.)
\bigskip
\hedd{Example.}
To compute the \BB/ transformation of a constant torsion curve, one needs to
solve the differential equation (\betaODE), which is equivalent to a Riccati
equation when one changes variables to $y=\tan(\beta/2)$.  When $\k$ is equal
to a constant $\k_0$,
the curve $\gamma$ is a helix, and the ODE is easy to solve.  For $C < \k_0$,
the solutions are periodic, and one can arrange that the period is commensurable
with the translational period of the helix, yielding a curve $\gnew$ that
is periodic up to translation (see Figure 1.1(b)).  For $C \ge \k_0$, the solutions
are asymptotically constant, and yield curves $\gnew$ that are asymptotic to
a helix (see Figure 1.1(c)).  For $C=\k_0$, the curvature of $\gnew$ is a rational
function of arclength:
$$\knew(s) = \k_0\dfrac{1-x^2}{1+x^2},\qquad x=\k_0 s.$$
\midinsert
\hskip -0.7truein
\BoxedEPSF{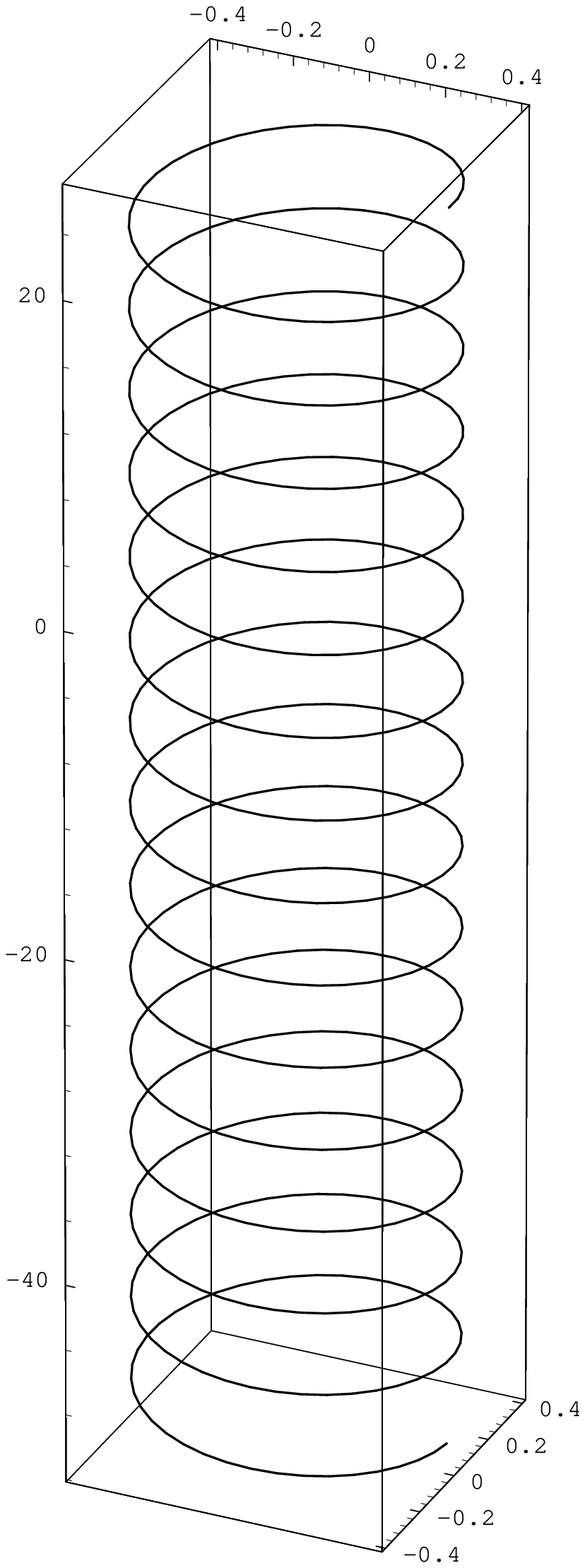 scaled 400}
\hskip -1.2truein
\BoxedEPSF{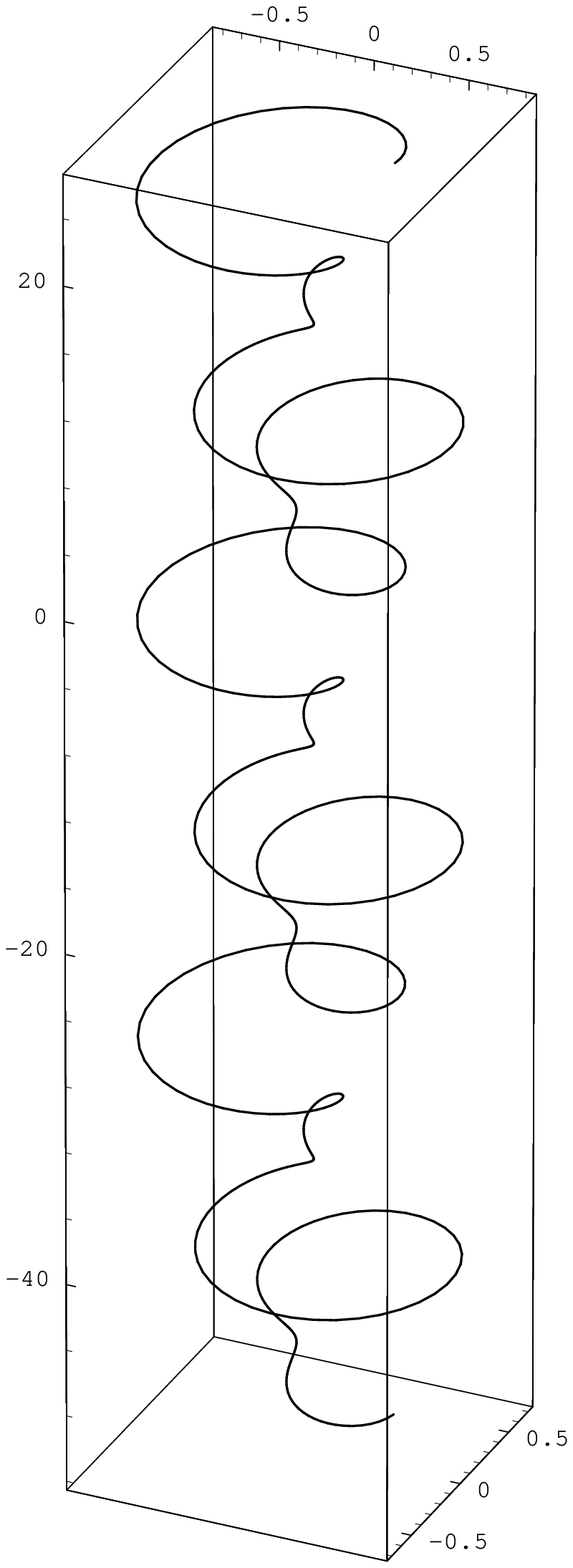 scaled 400}
\hskip -1.2truein
\BoxedEPSF{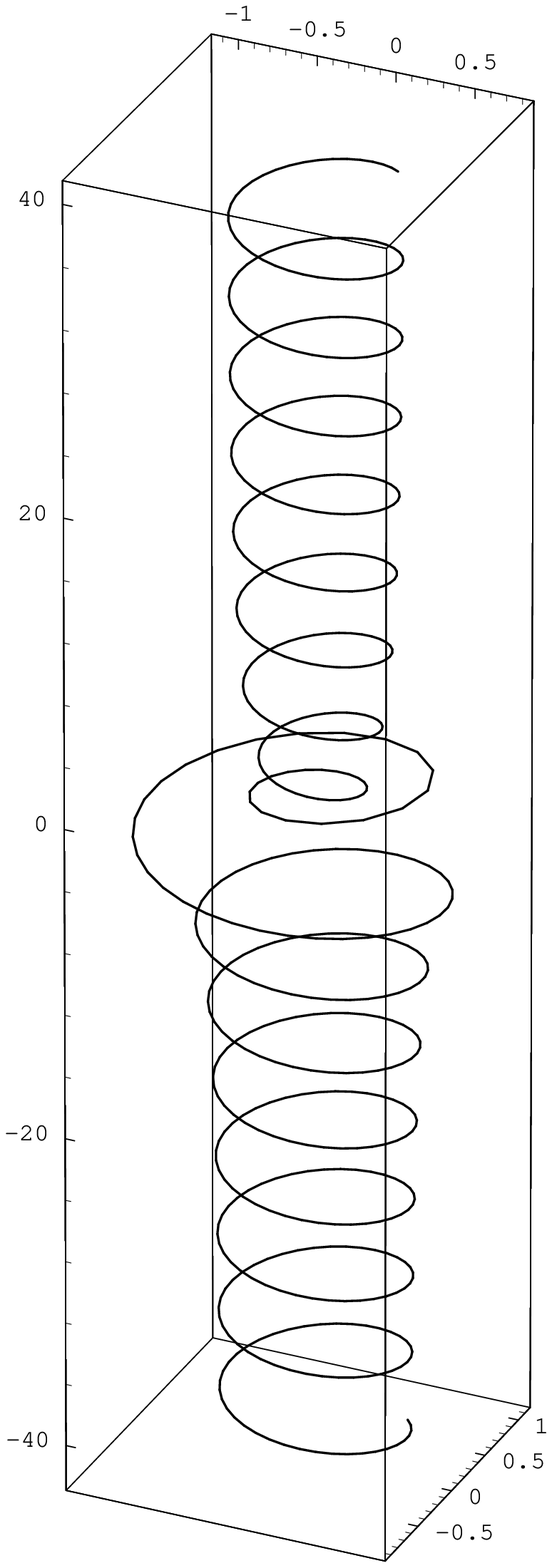 scaled 400}
\botcaption{Figure 1.1} A helix and its periodic and
 asymptotically helical \BB/ transformations \endcaption 
\endinsert
\bigskip
The \BB/ transformation for pseudospherical surfaces is, of course, equivalent to
the \BB/ transformation for the sine-Gordon equation.  If $s$ and $t$ are arclength
coordinates along the asymptotic lines, and $\theta$ is the angle between
the asymptotic lines, then the Codazzi equations imply $\theta_{xt} =
\sin\theta$.  In terms of $\varphi=\theta/2$, the \BB/ transformation for sine-Gordon
is
\define\barphi{\widetilde\varphi}
$$(\barphi -\varphi)_s = C\sin(\barphi+\varphi),\qquad (\barphi+\varphi)_t=
C^{-1}\sin(\barphi-\varphi).$$
If $T,N$ are chosen along an $s$-curve so as to agree with the orientation
$\di/\di t \wedge\di/\di s$ on the surface, then $\theta_s = \k(s)$.  Once
we make this identification,
then (\betaODE) follows if $\beta = -(\varphi+\barphi)$.

Along with applying the \BB/ transformation to constant torsion curves,
we can also apply the nonlinear superposition principle
for solutions of the sine-Gordon equation [Ro].
(For surfaces, this is Bianchi's theorem of permutability [E].)
Suppose $\beta_1 = \beta(s;\k(s),C_1)$ and $\beta_2 = \beta(s;\k(s),C_2)$
generate two \BB/ transformations, with corresponding curvature functions
$\k_1(s),\k_2(s)$ given by (\betak).  Then particular solutions
$\beta_{12}= \beta(s;\k_1(s),C_2)$ and $\beta_{21} =\beta(s;\k_2(s),C_1)$ to
successive \BB/ transformations are
given by
$$\beta_{12} - \beta_1 = \beta_{21} - \beta_2 = 2 \arctan\left[
\dfrac{C_1+C_2}{C_1 - C_2}\tan\left(\dfrac{\beta_1 - \beta_2}{2}\right)\right].
\tag\labeleq{\bianchi}$$

Now suppose that we can solve (\betaODE) for some initial curvature function and
all values of $C$, as we can for the helix.  Then in principle we can use
(\bianchi) to calculate arbitrarily many \BB/ transformations of the initial curve,
choosing constants $C_1, C_2, C_3, \ldots$.
(Recall also that once a single solution of a Riccati equation is known, all others
can be obtained by quadrature.)  The superposition principle (\bianchi)
also makes sense in
the limit as $C_2 \to C_1$, if we fix initial values for $\beta_1$ and $\beta_2$
independent of $C_1,C_2$.

\heading 1.2 The Frenet frame and the linear system \endheading
In this section we will see that, if we know the solution of the Frenet equations
for a given analytic curvature function $\kappa(s)$, 
as analytic functions of the constant torsion $\tau$, then we already
know the solution to the Riccati equation (\betaODE) that we have to solve to
compute the \BB/ transformation for the corresponding constant torsion curves.

For the given curvature function $\k(s)$, suppose $\Psi(s;\lambda)$ is the
fundamental solution
matrix for the linear system
$$\dfrac{d\psi}{ds} = \dfrac{1}{2}
\pmatrix \lambda &\k(s)\\-\k(s) &-\lambda \endpmatrix\psi.\tag\labeleq{\psisys}$$
(It follows from standard ODE theory that $\Psi$ is an analytic function of
$\lambda$.) When $\lambda = -i\tau$ for $\tau \in \R$, this system 
encodes the Frenet equations, in the following way.
If
$$\gamma(s)=i \Psi^{-1}\left. \dfrac{d\Psi}{d\lambda}\right|_{\lambda = -i\tau},
$$
then $\gamma$ is a curve in $su(2)$, which we will identify with $\R^3$ (see below).
One can verify that $\gamma$ satisfies
the Frenet equations for curvature $\k(s)$ and constant
torsion $\tau$, when the Frenet frame is given by the $su(2)$ matrices
$$\aligned
 T &= \dfrac{1}{2}\Psi^{-1}\pmatrix i & 0 \\ 0 & -i\endpmatrix\Psi \\
 N &= \dfrac{1}{2}\Psi^{-1}\pmatrix 0 & i \\ i & 0\endpmatrix\Psi \\
 B &= \dfrac{1}{2}\Psi^{-1}\pmatrix 0 & -1 \\ 1 & 0\endpmatrix\Psi\endaligned
\tag\labeleq{\symfr}$$

Now suppose that $\Psi(s;\lambda)$ were extended complex analytically from
$\lambda$ pure imaginary to all $\lambda \in \C$.  We observe then that
for $\lambda = C \in \R$, the coefficient matrix in (\psisys) is real, and
for any real vector solution $\psi$, the ratio of entries $y = -\psi_1/\psi_2$
satisfies
the Riccati equation
$$\dfrac{dy}{ds} = C y - \dfrac{\k}{2}(1+y^2).$$
Since this is the same as the equation satisfied by $\tan(\beta/2)$ under the
\BB/ transformation, we obtain a \BB/ transformation for the curve $\gamma$ using
$$\sin\beta = \dfrac{-\psi_1\psi_2}{\psi_1^2+\psi_2^2},\qquad
\cos\beta = \dfrac{\psi_2^2 - \psi_1^2}{\psi_1^2+\psi_2^2}.\tag\labeleq{\scbeta}$$

In practice, we may not know the fundamental matrix $\Psi$ explicitly, but
if we know the Frenet frame for the curve, normalized so that $(T,N,B)$ is
the identity matrix at $s=0$, and use the identification
$$(1,0,0) \to \dfrac{1}{2}\pmatrix i & 0 \\ 0 & -i\endpmatrix,\quad
(0,1,0) \to \dfrac{1}{2}\pmatrix 0 & i \\ i & 0\endpmatrix,\quad
(0,0,1) \to \dfrac{1}{2}\pmatrix 0 & -1 \\ 1 & 0\endpmatrix$$
of $\R^3$ with $su(2)$, then (\symfr) gives us values for all possible quadratics in the
entries of $\Psi$.  
Thus, for any linear combination $\psi$ of the columns of $\Psi$, we can express
the quadratics in (\scbeta) in terms of
 analytic extensions of those obtained using (\symfr).

\heading 1.3 Torus knots and closed elastic rods of constant torsion\endheading
In a recent paper [L-S4], Langer and Singer formulate and solve the Euler-Lagrange
equations associated to the {\it Kirchoff elastic rod}, by which is meant a curve that
is critical for some linear combination of $\int \k^2 ds$, $\int \tau ds$ and 
length.
Among the solutions they obtain are two
families of constant torsion curves.  For the family we will use, the curvature is given by an elliptic cosine, and the shape of the rod
is governed by the torsion, the maximum curvature $\k_0$, and the elliptic modulus $p$:
$$\k(s) = \k_0 \cn(x,p),\qquad\text{where } x = \dfrac{\k_0 s}{2p}.$$
(In what follows, we will take $\k_0=1$.)
Up to scale, then, the shape of the rod depends only on two parameters; we will
use the parameters $p$ and $\sigma = 2\tau$, and
the coordinate $x$ along the rod.  We will show how to choose the parameters to
obtain closed elastic rods of constant torsion.

Langer and Singer compute the position of the rod in a system of cylindrical
coordinates $(r,\theta,z)$ that are generated by the Killing fields associated
to the rod (see [L-S4], \S 5).  In
these coordinates, $r$ is already a $2K$-periodic function of $x$, where
$K=K(p)$ is the complete elliptic integral of the first kind:
$$\align r &= \dfrac{1}{\mu}\sqrt{\dfrac{1}{m} - \sn^2 x},\\
\intertext{where}
\mu &= \dfrac{1}{4}\sqrt{(p^{-2} - \sigma^2)^2 + 4\sigma^2},\\
m &= 16 \mu^2 / (p^{-2} + \sigma^2)^2.\endalign$$
\define\pmax{p_{\text{max}}}
Closure in $z$ imposes (cf. equation (26) in [L-S4])
$$\sigma^2 = \dfrac{1}{p^2}(2 \dfrac{E(p)}{K(p)} - 1),\tag\labeleq{\zclose}$$
where $E(p)$ is the complete elliptic integral of the second kind.  
(David Singer
has shown that constant torsion rods in the other family, whose curvature is
given elliptic $\dn$, are never closed in the $z$ coordinate.)
This implies $|\sigma| < p^{-1}$ and that $p$ has a maximum possible value
$\pmax \approx .9089085$ at which $\sigma^2$ approaches zero.  
(We will take $\sigma >0$,
giving positive torsion curves; curves obtained using $\sigma <0$ will differ
by a reflection.)  Using the remaining parameter $p$, we can attempt to make the change
in $\theta$ over a $2K$ period equal to a rational multiple of $2\pi$.  The
Euler-Lagrange equations imply that
$$\dfrac{d\theta}{dx} = -\dfrac{p\sigma}{\mu}\left( \lambda_1 +
\dfrac{(p^{-2}-1)(p^{-2}-\sigma^2)}{2(p^{-2} + \sigma^2)(1-m\sn^2 x)}\right),
\qquad\text{where } \lambda_1 =(\sigma^2 -p^{-2} +2)/4.$$
(This formula is due to Singer.)
So, the change in $\theta$ can be expressed in terms of complete elliptic integrals:
$$\align
-\dfrac{\Delta\theta}{2K} &= E(\xi,p') +
\left(\dfrac{E}{K}-1\right) F(\xi,p') + \dfrac{\lambda_1 p\sigma}{\mu}
 =: \Lambda(\sigma,p)
\tag\labeleq{\lamdef}\\
\intertext{where}
\xi &= \pi/2 - 2\arctan\left(\dfrac{2\sigma}{p^{-2} - \sigma^2 + 4\mu}\right)
\endalign$$
(cf. formula 413.01 in [B-F]).
\medskip
\hedd{Remark.}It will be desirable to have a formula for $\Lambda(p,\sigma)$ that does not
involve $\mu$ as denominator.  Using Legendre's relation, and the addition
formula and imaginary transformation for Jacobi's zeta function $Z(x)$ (see
[W-W], \S22.73), we obtain
$$\Lambda = -iZ\left(\sn^{-1}\alpha\right) + p^2\sigma\alpha,\qquad\text{ where\ }
\alpha = \dfrac{4\mu}{p(p^{-2}+\sigma^2)}.\
\tag\labeleq{\lambetter}$$
\medskip
The following result now follows directly from
the above formulas and those in [L-S4]:
\proclaim{Theorem \numbit}Given any relatively prime integers $m,n$ such
that $|m/n| <1/2$, there exists a smooth closed elastic rod of
constant torsion with the knot type of an $(m,n)$ torus knot.
\endproclaim
\demo{Proof}  We use continuity, noting that
the limiting values of $\Delta\theta/2\pi$ are zero
as $p\searrow 0$, and one half as $p \to \pmax$, to assert
that $\Delta\theta/2\pi=m/n$ for some modulus $p$.
Then one checks that,
over the course of one $2K$ period,
$r(x)$ and $z(x)$ describe a simple closed curve in the $rz$ plane.\qed\enddemo
\hedd{Remarks.}
It is interesting that this is exactly the same set of knot types as are available
among elastic curves (see [L-S1]), which of course are critical for a combination
of $\int \k^2 ds$ and length, and do not have constant torsion (unless they
are planar).
Since a computer-generated plot of $\Delta\theta$ shows it to be a monotone increasing function of $p$, we
expect that, similar to the situation for elastic curves, there is a unique 
constant torsion
elastic rod for each possible torus knot type.  (However, by arranging that
$\Delta\theta/2\pi=1/n$ for any $n\ge 3$, we can get many unknots.)
Some of these torus knots, as well as some unknots, are shown in Figure 1.2.
It is also apparent that, as $p \to \pmax$ and the torsion approaches zero,
the curve approaches a figure-eight elastic curve contained in a plane through the
$z$-axis.
\midinsert
\hskip 0.2in
\TrimTop{2.6in}\TrimBottom{2.6in}
\BoxedEPSF{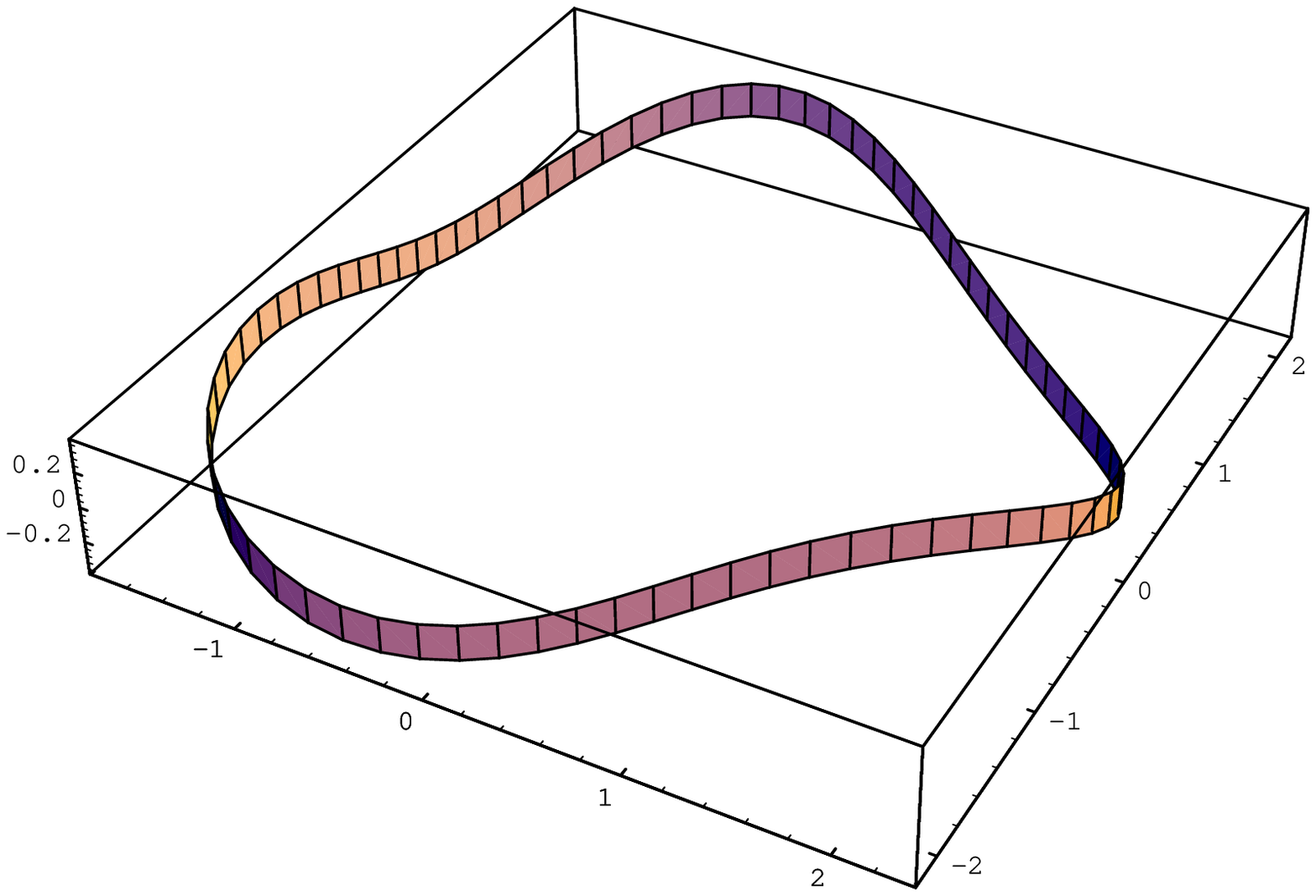 scaled 300}
\hskip 0.2in
\TrimTop{2.6in}\TrimBottom{2.6in}
\BoxedEPSF{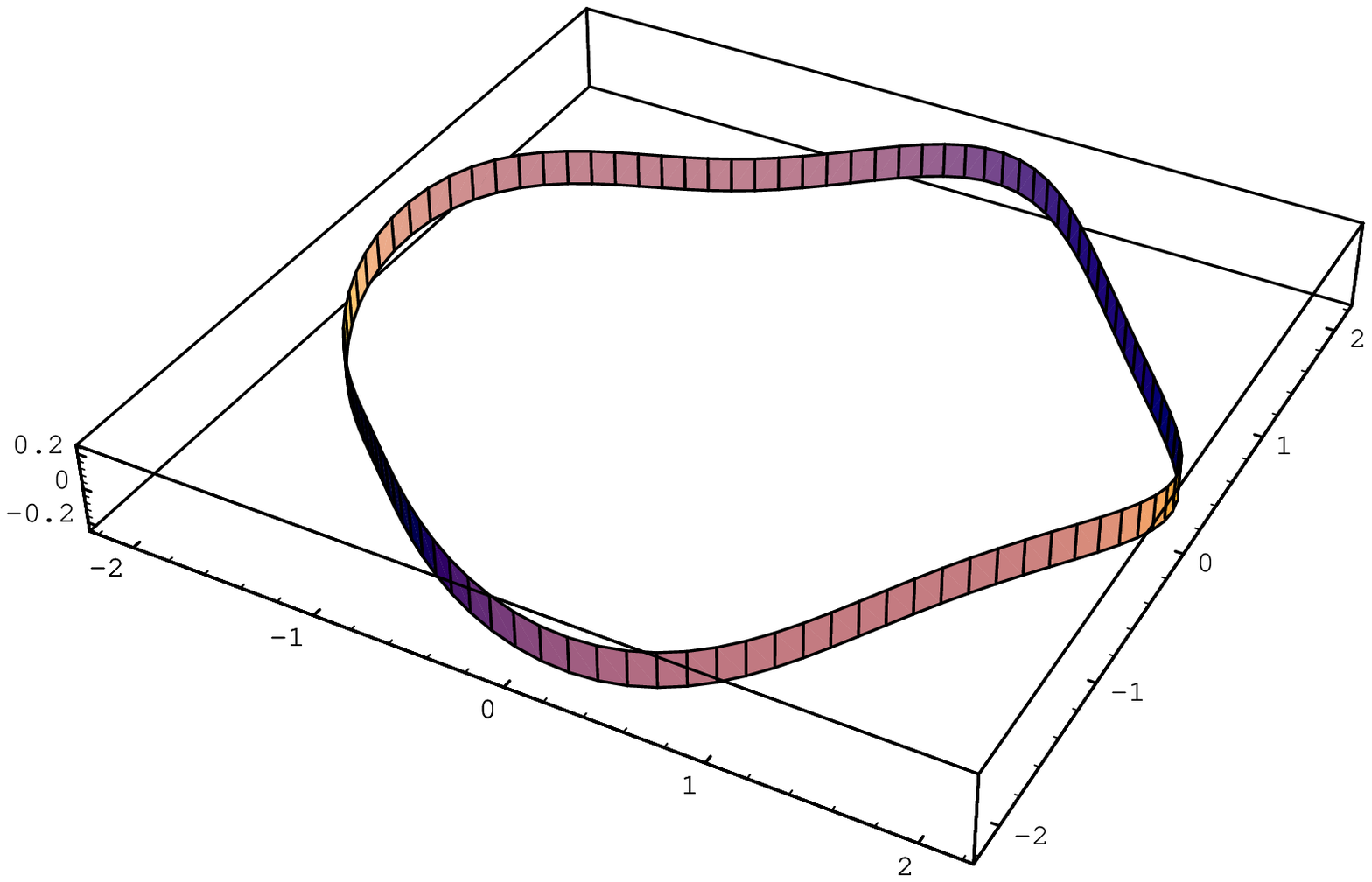 scaled 300}
\bigskip
\hskip 0.2in
\TrimTop{2.6in}\TrimBottom{2.6in}
\BoxedEPSF{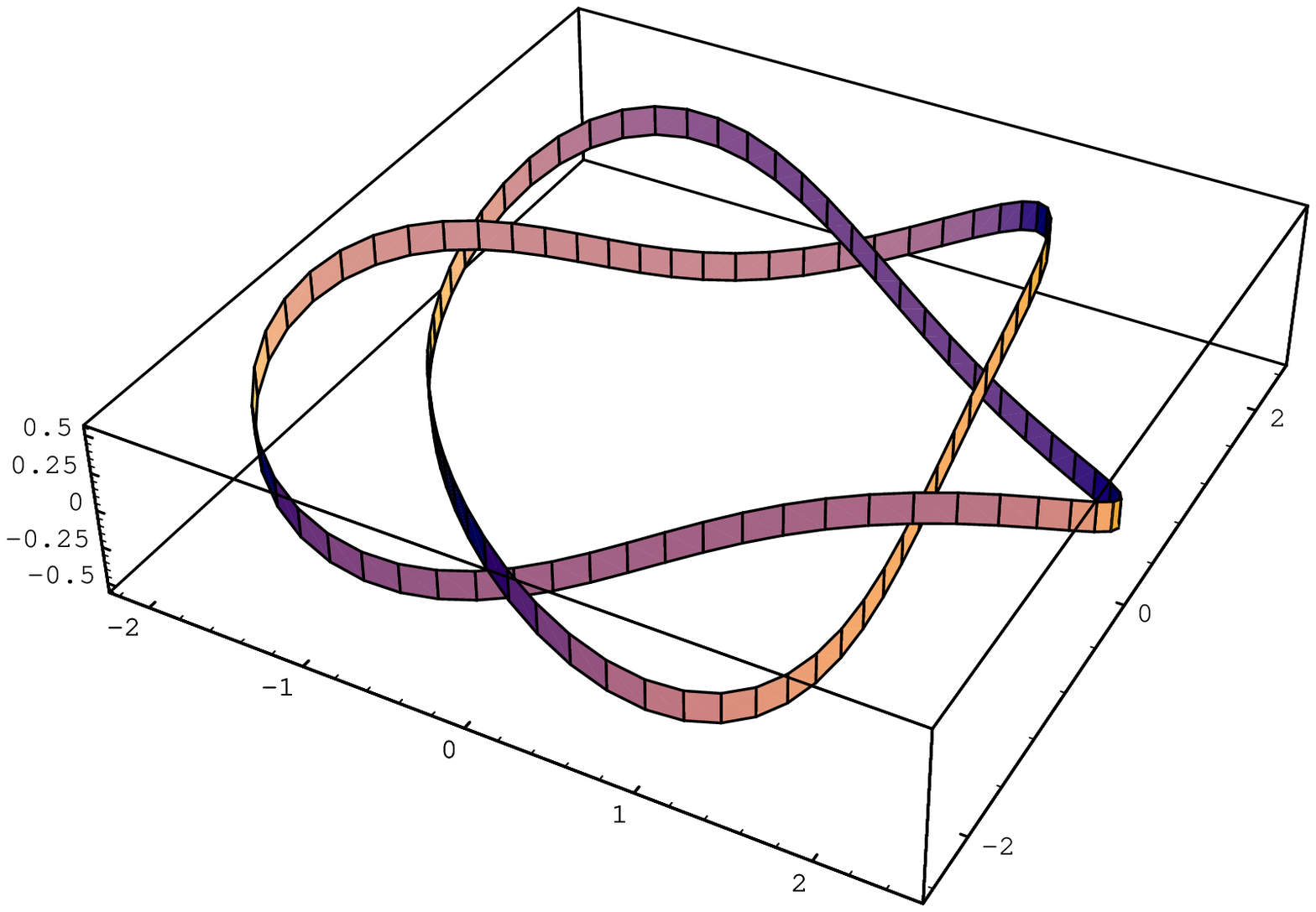 scaled 300}
\hskip 0.2in
\TrimTop{2.6in}\TrimBottom{2.6in}
\BoxedEPSF{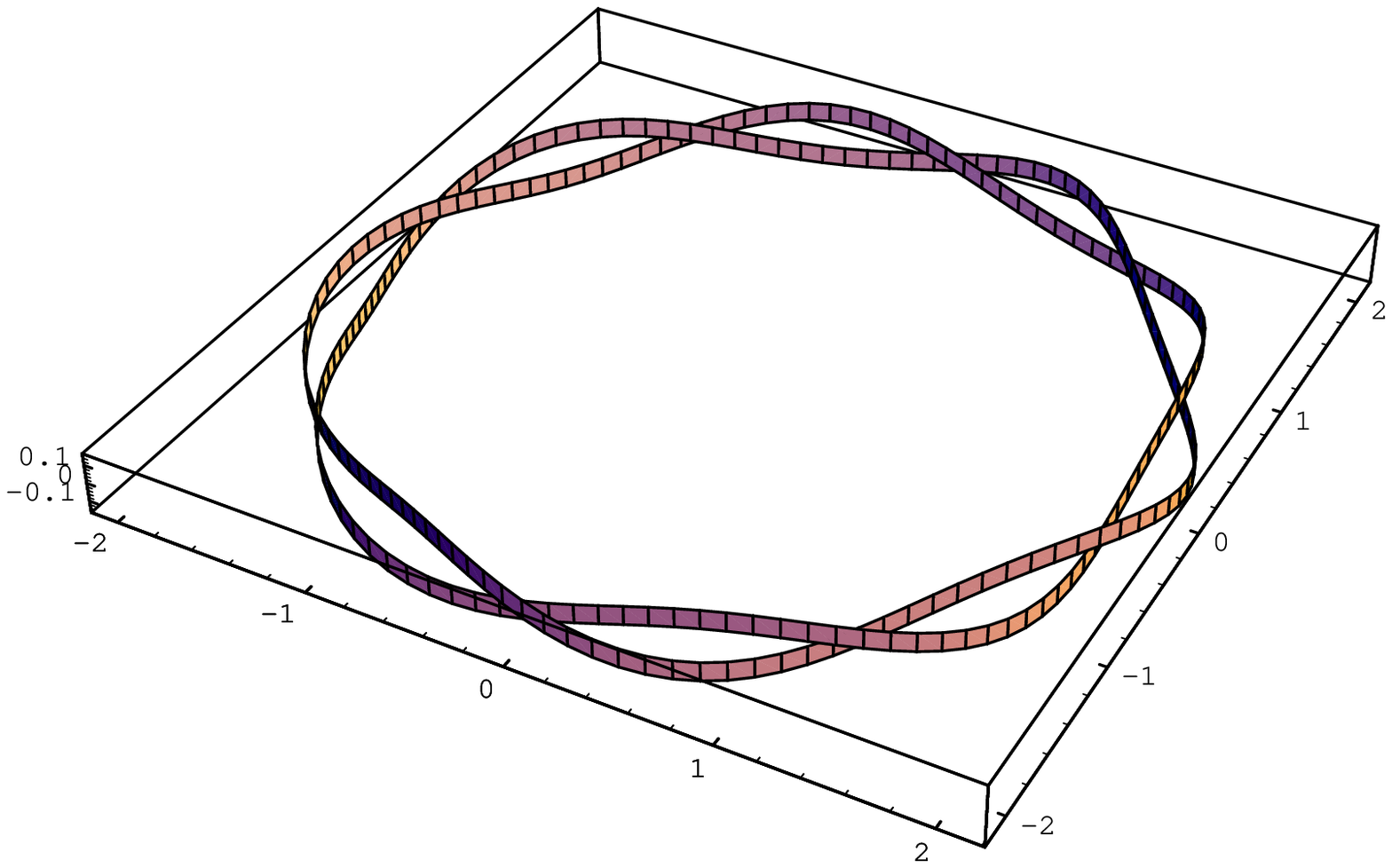 scaled 300}
\bigskip
\hskip 0.2in
\TrimTop{2.6in}\TrimBottom{2.6in}
\BoxedEPSF{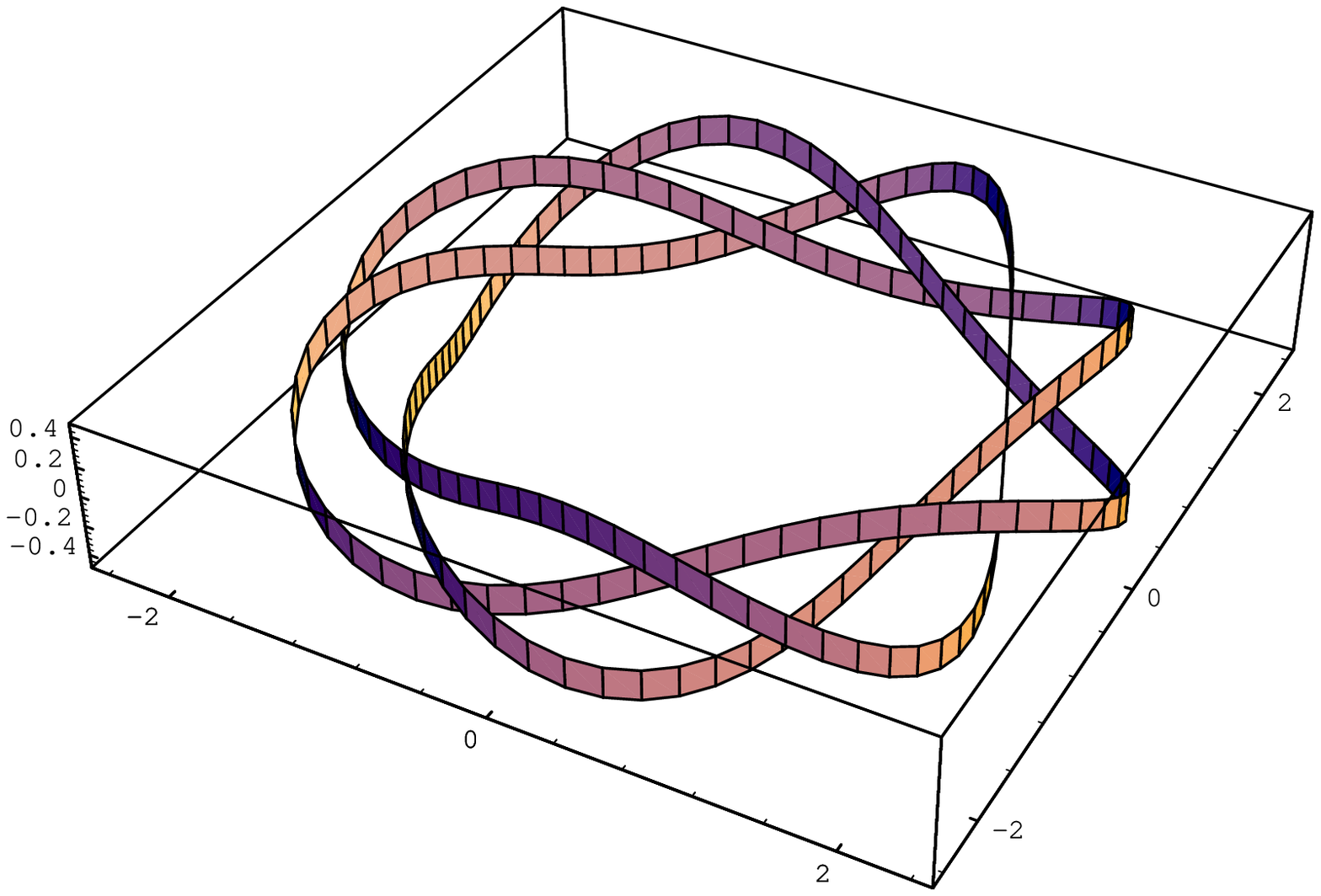 scaled 300}
\hskip 0.2in
\TrimTop{2.6in}\TrimBottom{2.6in}
\BoxedEPSF{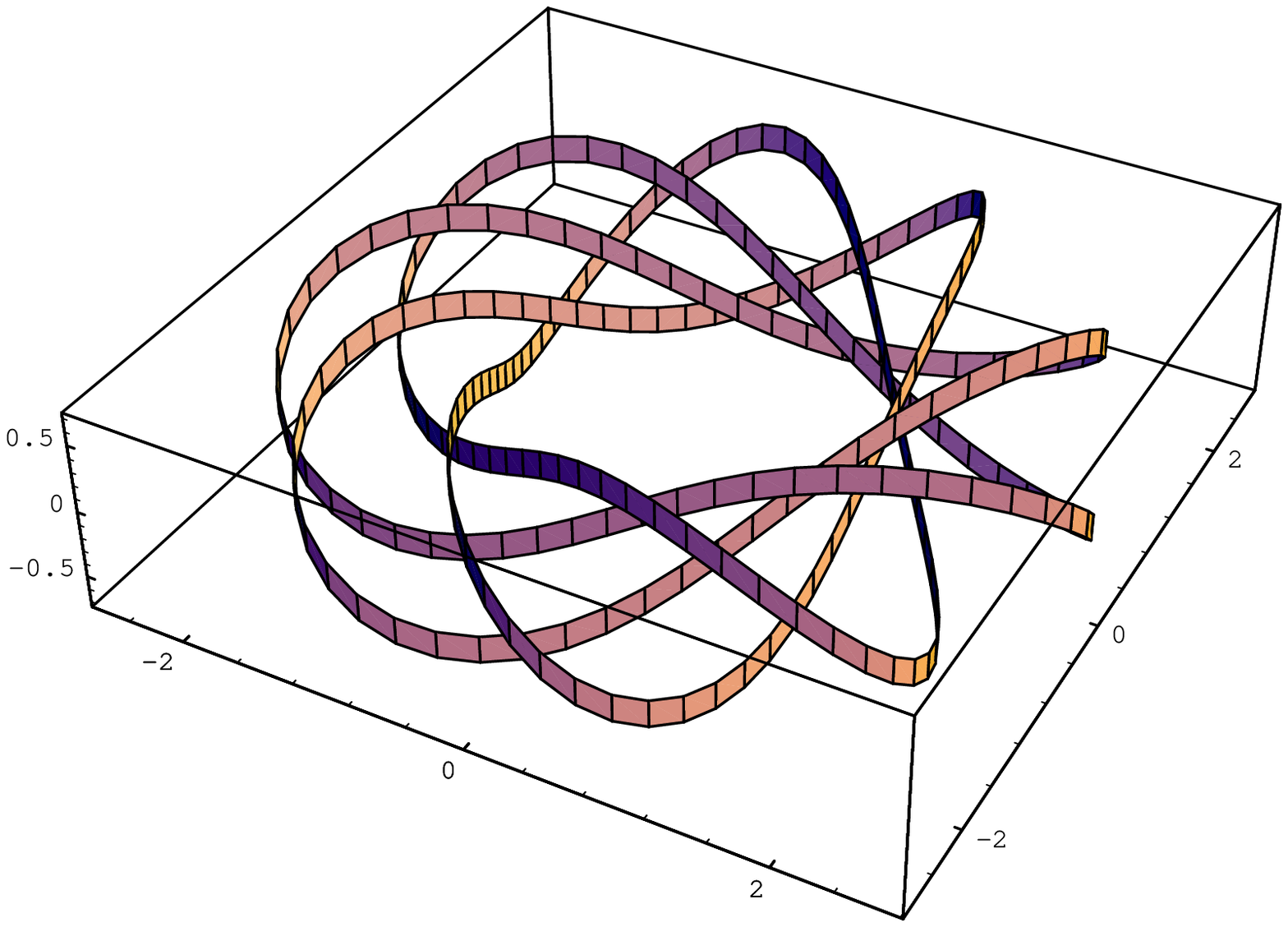 scaled 300}
\botcaption{Figure 1.2} Some closed elastic rods of constant
torsion.  (The curves here have been thickened by a vertical ribbon.)
\endcaption
\endinsert
\bigskip
For the purpose of computing their \BB/ transformations, we need to compute
Frenet frames for these elastic rods.  Using the formulas in [L-S4], one can 
express the Frenet frame in a cylindrical coordinate basis:
$$\pmatrix T\\N\\B\endpmatrix =
\pmatrix -\dfrac{\cn x\sn x\dn x}{2p\mu^2r^2} &
                               -\dfrac{\sigma(p^{-2}+\sigma^2-4\lambda_1 \sn^2 x)}{8\mu^3 r^2} &
         \dfrac{\cn^2 x -2\lambda_1}{2\mu} \\
                       -\dfrac{\cn x(p^{-2}+\sigma^2 - 2\sn^2 x)}{4\mu^2 r^2} &
                              \dfrac{\lambda_1 \sigma \sn x\dn x}{2p\mu^3 r^2}&
         -\dfrac{\sn x\dn x}{2 p \mu} \\
         \dfrac{\sigma \sn x \dn x}{2 p \mu^2 r^2}&
         -\dfrac{(p^{-4} - \sigma^4)\cn x}{16\mu^3 r^2}&
         -\dfrac{\sigma \cn x}{2\mu} \endpmatrix
\pmatrix r\di_r\\ \di_\theta \\ \di_z \endpmatrix
$$
To convert to Cartesian coordinates, we need the cylindrical coordinates
$r,\theta$ as functions along the elastic rod.  Using 
formulas that express incomplete elliptic integrals of the third kind in terms
of theta functions (cf. formula 434.01 in [B-F]), we get
$$\aligned r\cos\theta &= \sqrt{\dfrac{2Kp'}{\pi p}}\
\dfrac{e^{-i\Lambda x}\Theta_1(x-i\widehat F) + e^{i\Lambda x}\Theta_1(x+i\widehat F)}
{2\mu\Theta(x)H_1(i\widehat F)}\\
r\sin\theta &= \sqrt{\dfrac{2Kp'}{\pi p}}\
\dfrac{e^{-i\Lambda x}\Theta_1(x-i\widehat F) - e^{i\Lambda x}\Theta_1(x+i\widehat F)}
{2\mu i\Theta(x)H_1(i\widehat F)}\endaligned\tag\labeleq{\rodxy}$$
where $\widehat F = F(\xi, p')$ 
and $\Lambda=\Lambda(\sigma,p)$ is defined in (\lamdef).

Now let $H$ be the matrix whose rows are $T,N,B$.  Let $H_0 = H\restr_{x=0}$, and let $G(x,\sigma) = H_0^{-1} H$; $G$ is the normalized Frenet
frame.
The above formulas can be used to express the entries of $G$ as complex analytic functions of
$\sigma$ in a disk about the origin in the $\sigma$-plane.  (One can also
desingularize the formulas near the points $\sigma = \pm i \pm p'/p$, where $\mu=0$.)
In particular, from $G(x,2iC)$ we may, using the formulas (\scbeta) in the
previous section, construct the \BB/ transformation of an elastic rod of constant
torsion. However, the new curve is not necessarily closed; in subsequent sections
we will see how to obtain closure and how severely closure restricts the shape
of the new curve.
\medskip
\hedd{Remarks.}  

1. From the formulas in [L-S4] one can also obtain 
$$z(x) = \dfrac{1}{\mu p} \dfrac{\Theta'(x)}{\Theta(x)}.\tag\labeleq{\rodz}$$
This bears a striking resemblance to the formula obtained by Mumford [M] for
planar elastic curves in terms of theta functions.  In fact, one can check that
as $p \to \pmax$ and our elastic rod becomes planar, (\rodxy) and (\rodz) 
agree with Mumford's formula in the limit.

2. One can derive formulas quite similar to  (\rodxy) and (\rodz) for elastic rods 
in general.

\heading 1.4 Closure of the \BB/ transformation \endheading
Suppose we have generated a \BB/ transformation of a
constant torsion curve $\gamma$,
by obtaining the fundamental matrix $\Psi$ for the linear system (\psisys), selecting
a non-zero vector solution
$$\psi = \Psi \pmatrix c_1 \\c_2 \endpmatrix \qquad c_1, c_2 \in \R,$$
and performing the transformation using (\scbeta) and (\curveback).  Suppose that
$\gamma$ is closed and has length $L$.  Then the \BB/ transform $\gnew$ will
be closed if the vector $V=\cos\beta T+\sin\beta N$ has period $kL$ for $k$ a positive
integer.

\proclaim{Proposition \numbit} Suppose $V$ has period $kL$.  If $k$ is even or
$\gamma$ is even, then $(c_1, c_2)$ is an eigenvector of $\Psi\restr_{s=kL}$.
If $k$ is odd and $\gamma$ is odd, then $(c_1,c_2)$ is an eigenvector of
$\left( \smallmatrix 1&0\\0& -1\endsmallmatrix \right) \Psi\restr_{s=kL}$.
\endproclaim

When $\gamma$ is an elastic rod of constant torsion, we may use the entries of
the normalized Frenet frame $G(x,\sigma)$, together with the formulae (\symfr), to obtain
all quadratics in the entries of $\Psi$.  In particular, to obtain the eigenvectors
required in the above proposition, we may use the squares of the relevant
matrices.

Suppose that the elastic rod closes up after $n$ periods of length $2K$ in the
$x$ parameter.  (The rod is even or odd exactly as $n$ is even or odd.)  If $kn$
is even, then we calculate that
$$\Psi^2\restr_{x=2knK} = \cos(x\Lambda) \pmatrix 1&0\\0&1 \endpmatrix
- \sin(x \Lambda) \pmatrix  i\dfrac{1-2\lambda_1}{2\mu} & \dfrac{\sigma}{2\mu} \\  -\dfrac{\sigma}{2\mu}& -i\dfrac{1-2\lambda_1}{2\mu}
\endpmatrix$$
where $\Lambda$ is given by (\lamdef).
Note that, when we want to calculate \BB/ transformations, we set 
$\sigma = 2\tau = 2iC$ for $C\in \R$, $C\ne 0$, and then
$\Lambda$ is
purely imaginary and $\Psi^2$ is real.  Moreover, the eigenvalues
$\cos(x\Lambda) \pm i\sin(x\Lambda)$ of this matrix are always distinct.
If $kn$ is odd, then
$$\left(\left( \smallmatrix 1&0\\0& -1\endsmallmatrix \right)
\Psi\right)^2\restr_{x=2knK} =  -\cos(x\Lambda) \pmatrix 1&0\\0&1 \endpmatrix
+ \sin(x \Lambda) \pmatrix  i\dfrac{1-2\lambda_1}{2\mu} & \dfrac{\sigma}{2\mu} \\  -\dfrac{\sigma}{2\mu}& -i\dfrac{1-2\lambda_1}{2\mu}
\endpmatrix$$
---i.e., the same matrix as above, up to a minus sign.  Thus, in either case
the eigenvalues are distinct.  The corresponding eigenvectors are
$$v_+ = \pmatrix -\dfrac{\sigma}{2\mu}\\i\left( 1 + \dfrac{1-2\lambda_1}{2\mu}\right)
\endpmatrix,\qquad
v_- = \pmatrix i\left( 1 + \dfrac{1-2\lambda_1}{2\mu}\right)\\
-\dfrac{\sigma}{2\mu}\endpmatrix.\tag\labeleq{\eigens}$$
Since these eigenvectors are independent of $k$, we conclude that
{\it when $\gamma$ is a closed elastic rod of constant torsion,
$\gnew$ closes up after $k$ circuits around $\gamma$ if and only if it
closes up after one circuit.}

\heading 1.5 Elastic rods again \endheading
In this section we will see that, for any value of the \BB/ parameter $C$,
a closed \BB/ transform of a closed elastic rod of constant torsion is congruent
to that rod.  Since the two curves have the same torsion, it will suffice to
show that they have the same curvature function, up to a phase shift.

Let $q = 2C$ and $\sigma = iq$ for $q \in \R$.  Then by (\scbeta)
the initial values for $\beta$
corresponding to eigenvectors $v_+$ and $v_-$ given above are
$$\sin \beta = \dfrac{2q}{p^{-2} + q^2},\qquad \cos\beta =
\pm\dfrac{\sqrt{(p^{-2}+q^2)^2-4q^2}}{p^{-2} + q^2}.\tag\labeleq{\binit}$$
Since $\k = \cn x$, the differential equation for $\beta$ is
$$d\beta/dx = p(q\sin \beta - 2\cn x).\tag\labeleq{\dbdx}$$

If $\knew$ were $\cn(x-a)$ then (\betak) would imply that
$$q\sin \beta = \cn x - \cn(x-a)\tag\labeleq{\qsb}$$
and
$$d\beta/dx = -p(\cn x + \cn(x-a)).$$
Then by integration,
$$\align \beta &= -[\sin^{-1}(p\,\sn x) +\sin^{-1}(p\,\sn(x-a))] + C_1
\tag\labeleq{\betasol}\\
\intertext{for some constant $C_1$, giving}
\sin(\beta-C_1) &= -p[\sn x \dn(x-a) + \dn x\sn(x-a)]\\
\cos(\beta-C_1) &= \dn x \dn(x-a) -p^2 \sn x \sn(x-a).
\endalign$$
We can find the constants $a$ and $C_1$ by setting $x=0$: (\qsb) gives
$$\cn a = (p^{-2} - q^2)/(p^{-2} + q^2)$$ and so
$$\sn a = \pm \dfrac{2p^{-1}q}{p^{-2} + q^2},\tag\labeleq{\sna}$$
while (\betasol) gives $\sin(\beta-C_1) = p \sn a$ when $x=0$.  Comparing with
(\binit) shows that we should use $C_1=0$ or $C_1=\pi$; since $\cos(\beta-C_1)$
is always positive in (\betasol), we use $C_1=0$ and the plus sign
in (\sna) for $v_+$, $C_1=\pi$ and the minus sign in (\sna) for $v_-$.

It remains only to be seen that (\betasol) gives a solution to (\dbdx). (We will
show this for the initial values associated to $v_+$, using $C_1=0$ and
the plus sign in (\sna); the $v_-$ case is similar.)  Let $b=a/2$.  Then from
(\betasol),
$$d\beta/dx =-2p\cn x -2p\left(
\dfrac{\sn b\,\dn b\, \sn(x-b)\dn(x-b)}{1-p^2 \sn^2 b\,\sn^2(x-b)}\right)$$
(cf. 123.03 in [B-F]).  Since
$$pq \sin \beta = -2p^2 q\dfrac{\sn(x-b)\dn(x-b)\cn b}{1-p^2 \sn^2 b\,\sn^2(x-b)},$$
we will be done if we can show that $pq \cn b = \sn b\dn b$.  This now follows
from the double angle formula
$$\dfrac{\sn a}{1+\cn a} = \dfrac{\sn b\, \dn b}{\cn b}.$$
Now, by construction, $\knew = \cn(x-a)$ and hence $\gnew$ is congruent to
the elastic rod $\gamma$.

\heading 1.6 Linking numbers \endheading
In this section we use the explicit solution of (\betaODE) obtained for an
elastic rod $\gamma$ in the previous section, to compute the linking number
of $\gamma$ and its twin $\gnew$ in terms of the {\it self-linking} of
$\gamma$ (see [Po1] for discussion of self-linking).

\proclaim{Theorem \numbit} If $|C|$ is sufficiently small,
$$Lk(\gamma,\gnew) = SL(\gamma) - n/2\tag\labeleq{\linker}$$
where $n$ is the number of $2K$ periods of $\gamma$.
\endproclaim
The proof will be an application of White's formula (see [W], [Po2]).
Before proving this result, we should remind the reader that the self-linking
number $SL(\gamma)$ can be thought of as the linking number of $\gamma$ with
the curve traced out by the endpoint of the Frenet normal $N$ (assuming
$\gamma$ is even), suitably scaled so that the vector never intersects
other points of $\gamma$.  As such, $SL(\gamma)$ is the sum of
a contribution from the twisting
of $N$ about $\gamma$ and a contribution from the twisting of $\gamma$ about
itself---a measure of the knottedness of $\gamma$.  In (\linker), the former contribution
is cancelled out by the $-n/2$.  In particular, when $\gamma$
is unknotted (as it is when, for example, we arrange that $\Delta\theta = \pi/2$),
then $Lk(\gamma,\gnew)=0$.

\define\CC{\gamma}
\define\CCV{\gamma_{\lower.5em\hbox{V}}}
\demo{Proof}Let us recall White's Formula.  Let $\CC$ be a smooth closed space curve,
oriented by unit tangent vector $T$, and $V$ a unit normal vector along $\CC$.
Let $\CCV$ be the curve traced out by the endpoints of $\delta V$, with the
obvious orientation, where $\delta >0$ is chosen small enough that the ribbon
spanned by $\delta V$, with boundary $\CC \bigcup \CCV$, is embedded.  Then
$$Lk(\CC,\CCV) = Wr(\CC) + {1\over 2\pi}\int (T\times V)\cdot dV,
\tag\labeleq{\whiteform}$$
where the {\it writhe} is given by the integral
$$Wr(\CC) = \underset{\CC \times \CC}\to{\int \int} e^* dS,$$
wherein for $(x,y) \in \CC\times\CC$, $x\ne y$, $e$ is the unit vector from
$x$ to $y$, and $dS$ is the element of area on the unit sphere.  Before White,
Pohl [Po1]
obtained a special case of this formula,
$$SL(\CC) = Wr(\CC) + {1\over 2\pi}\int \tau ds,\tag\labeleq{\pohlform}$$
which allows us to express the writhe in terms of the self-linking of $\CC$.
(While Pohl did not define the self-linking number for curves with inflection
points, we will take (\pohlform) as defining $SL(\CC)$ for such curves.  If $\CC$ is an
even curve, this $SL(\CC)$ is indeed the linking number of $\CC$ and
$\CC_N$---for either choice of $N$---while if $\CC$ is odd, $SL(\CC)$ is a
half-integer.)

We next remark that it is easy to extend (\whiteform) to the case where $V$ is
a transverse vector field along $\CC$:
$$Lk(\CC,\CCV) = Wr(\CC) + {1\over 2\pi}\int (T\times U)\cdot dU,
\tag\labeleq{\whiteformU}$$
where $U$ is the unit vector in the direction of the orthogonal projection of $V$
into the normal plane along $\CC$.  However, for the \BB/ transformation we have
$$V= \cos\beta\ T + \sin\beta\ N,$$
which fails to be transverse to $\gamma$ whenever $\sin\beta =0$.  In order to
use (\whiteformU) to calculate the linking number, we will take a small perturbation
$\widetilde V$ of $V$ that is transverse to $\gamma$, also ensuring that
the ribbon spanned by $\widetilde V$ is embedded.  Then, since, as shown by
the Gauss integral formula, the linking number depends continously on the two curves,
our calculation of $Lk(\CC,\CC_{\widetilde V})$ will give $Lk(\CC,\CCV)$ also.

We may assume that $\kappa \ne 0$ and $\tau \ne 0$ in the vicinity of any point
where $V$ is tangent to $\gamma$.  A local calculation then shows that, when
$V$ is a positive multiple of $T$ at the point of tangency, the perturbation
$$\widetilde V = V -\epsilon \k\tau B,\qquad \epsilon > 0
\tag\labeleq{\perturb}$$
makes the ribbon spanned by $\delta \widetilde V$ embedded, for $\delta$ and
$\epsilon$ sufficiently small.  (It is interesting to note that if $\epsilon$
is negative in (\perturb), then no matter how small $\delta$ is chosen, the
ribbon fails to be embedded near the point of tangency.  The reader may visualize
this by imagining a momentarily tangent vector field $V$ along a right-handed
helix, with $\kappa >0$ and $\tau >0$: the vector field must be perturbed down,
not up, in order to prevent the perturbed vector field from intersecting the
helix.)   Now (\whiteformU) gives
$$Lk(\CC,\CCV) = SL(\CC) +
{1\over 2\pi} \int d \arctan\left(\dfrac{\sin \beta}{\epsilon \k\tau}\right)$$

When $\gamma$ is a constant torsion elastic rod, and the \BB/ transform is
constructed using eigenvector $v_+$,
$$\sin \beta = \dfrac{\cn x - \cn(x-a)}{q},\qquad \sn a = \dfrac{2pq}{1+p^2q^2}$$
and $\cos\beta >0$ always, showing that our assumptions hold.  Furthermore,
$\sin \beta/\k$ goes from $+\infty$ to $-\infty$ between inflection points, and
hence
$$Lk(\CC,\CCV) = SL(\CC) -n/2.$$
The same result holds when we use eigenvector $v_-$, but since $\cos \beta <0$
we have to use $\epsilon <0$ in (\perturb); however,
$\sn a = -\dfrac{2pq}{1+p^2q^2}$ makes up for the change in sign.
\qed\enddemo

\heading 1.7 \BB/ transformations and the LIE hierarchy \endheading
The constant torsion elastic rods are among the ``soliton'' curves for the 
Localized Induction Equation (LIE), and as such are critical for some 
linear combination of the
integral invariants associated with the LIE.  
In this section, we will formulate
a precise conjecture as to how the \BB/ tranformation takes solitons to 
solitons, and in particular how the coefficients in the linear combination of 
invariants change. Ercolani, Forest and Mclaughlin [E-F-M] use arguments
from algebraic geometry to address the analogous question for n-soliton
solutions of the sine-Gordon equation. They find that in general a 
\BB/ transformation (\betak ) takes n-solitons to (n+1)-solitons; 
what we describe below is a more geometric argument supporting the same 
conclusion.  In the
particular case when the initial condition of the Riccati equation
selects one eigenvector of the transfer matrix of the linear system,
producing closed curves as we discussed above, then [E-F-M] find that
the \BB/ formula takes an n-soliton to another n-soliton; this of course
agrees with our results on closed \BB/ transformations of elastic rods.

We begin with the observation, due to Langer and Perline [Pe], that along a
curve of constant torsion $\tau$ the linear combinations of the 
LIE vector fields $\{ X_k \}_{k=0}^{\infty}$ defined by
$$B_n = \sum_{k=0}^{2n} \binom {2n}k (-\tau)^{2n-k} X_{k+1},
\qquad n\ge 0$$
are purely binormal and do not involve $\tau$.  For example,
$$B_0 = X_1 = \k B,\qquad B_1 = X_3 - 2\tau X_2 + \tau^2 X_1 = -(\k''+k^3/2)B.$$
A {\it planar-like n-soliton} [Pe] is a curve of constant torsion along which
some constant-coefficient linear combination $\sum_{k=0}^n a_k B_k$
 vanishes.
For example, a helix is a planar-like 1-soliton.
Because the LIE vector fields $X_k$ are Hamiltonian vector fields associated
with certain densities involving $\k,\tau$ and their derivatives (see [L-P1] for
details), these solitons are critical curves for a natural variational problem.
For example, since $X_1,X_2,X_3$ are Hamiltonian vector fields for $\int ds$, 
$\int \tau ds$, 
$\int {1\over 2}k^2 ds$, an elastic rod of constant torsion is also a
planar-like 1-soliton.

It is easy to show, using the LIE recursion operator [L-P2], that if along
a constant torsion curve $\gamma$ we define
$$U_n = \sum_{k=0}^{2n} \binom {2n}k (-\tau)^{2n-k} X_k,
\tag\labeleq{\Udef}$$
then $\sum_{k=0}^n a_k B_k=0$ if and only if
$$J=\sum_{k=0}^n a_k U_k\tag\labeleq{\Jdef}$$
is a constant vector.  We call $J$ a {\it Killing field along $\gamma$} because
it is the restriction of a vector field generating a symmetry of Euclidean
space, in this case translation.  Combinations of LIE vector fields that are
Killing fields play an important role in showing that variational problems
associated to combinations of the LIE functionals yield finite-dimensional
completely integrable Hamiltonian systems (see [L-S2], [L-S3], [L-S4]).  So, it
is important to see how Killing fields along planar-like $n$-solitons interact
with the \BB/ transformation.

\proclaim{Conjecture \numbit} Let $\gamma$ be a constant torsion $n$-soliton, and let
$\gnew$ be a \BB/ transformation of $\gamma$.  Let $\overline U_k$ denote the
combinations of LIE vectors defined by (\Udef), computed along $\gnew$.  Then
if $J$ is the constant vector defined by (\Jdef) along $\gamma$,
$$J = \dfrac{1}{\tau^2 +C^2} \left[ \sum_{k=0}^n a_k \overline U_{k+1}
+b \overline U_0 + C^2\sum_{k=0}^n a_k \overline U_k\right]$$
where $C$ is the parameter in (\betaODE), (\curveback) and $b$ is a constant
depending only on $\gamma$.
\endproclaim

For any given $n$, one can mechanically verify this conjecture using the 
formulae for the vector fields $U_k$, $\overline U_k$, the equation 
(\betak) for the curvature
of $\gnew$, the differential equation (\betaODE), and the ODE for the
curvature of $\gamma$ implied by $\sum_{k=0}^n a_k B_k=0$.  (The constant $b$
arises as a first integral of this equation.)

We have verified the conjecture for $n=0,1,2$; hence we know that the \BB/
transform takes a planar-like 1-soliton to a planar-like 2-soliton, and takes
a planar-like 2-soliton to a planar-like 3-soliton.
However, the results of previous sections show that a {\it closed} transform 
of a closed 1-soliton is still a 1-soliton; an intuitive explanation for
this can be given for small values of the arbitrary constant
$C$ in the \BB/ formula. When $C \rightarrow 0$, since $\sin \beta$ is
bounded, the Riccati equation becomes approximately $\displaystyle
\frac{d\beta}{ds} \simeq -\k$. Therefore $\beta(s)\simeq -\theta(s) =
- \int^s \k(s')ds'$ and, to first order,
$$
\tilde{\gamma} = \gamma + \frac{2C}{\tau^2} \left[ \cos\theta T - 
\sin \theta N \right]. $$
In other words, the \BB/ transformation flows along the vector field
$$
\gamma_t \simeq W = \cos \theta T - \sin \theta N.
$$
This vector field $W$, discussed
in [L-P2], is known as the {\it trigonometric vector field}, and the associated
evolution for $\theta = \int^s \k(s')ds'$ is described by the
sine-Gordon equation, which {\it preserves} soliton type. Moreover,
for periodic solutions of the sine-Gordon equation, the limiting \BB/
formula described above transforms closed n-soliton curves into
closed n-soliton curves. 

In the second part of this paper, we will
perform iterated \BB/ transformations on the 1-solitons to obtain new
closed curves that are multi-solitons for the LIE hierarchy.

\heading II. The double \BB/ transformation \endheading
\newsec{2}
In Part I we used B\"{a}cklund transformations to produce new constant
torsion elastic rods from given ones. Moreover, we verified that
imposing the closure condition on the transformed curves restricts
them to the same congruence class as the original elastic rod.

More exotic and interesting constant torsion curves are produced
by successive iterations of the B\"{a}cklund transformation
starting with a given elastic rod of constant torsion.

\heading 2.1 A gauge formula for the \BB/ transformation \endheading
We present an alternative derivation of the \BB/ formula based on
a gauge transformation for the associated linear system, as described in
[E-F-M]. 
This approach is most convenient for constructing iterated \BB/ 
transformations, since  it reduces each step of the iteration to
a purely algebraic computation. We will show that the Riccati equation
(equivalently, the linear system) needs to be solved only once at the
initial step of the iteration.

We begin with the following form of the associated linear system
$$\dfrac{d\phi}{ds} = \dfrac{1}{2}
\pmatrix i\k(s) &\lambda \\ \lambda &-i\k(s) \endpmatrix\phi,
\tag\labeleq{\phisys}$$
which simplifies the formulation of the \BB/ transformation.
We observe that equation (\phisys ) is transformed into 
the linear system (\psisys ) by setting 
$\phi =A\psi$, where
$A = \dfrac{1}{2}\left( \smallmatrix 1-i & -1-i \\1-i & 1+i \endsmallmatrix \right)$.

Using the fundamental solution matrix $\Phi$ of (\phisys ), we
identify the curve of constant torsion $\tau$  with the
$su(2)$ matrix
$$
\gamma(s)=i \Phi^{-1}\left. \dfrac{d\Phi}{d\lambda}\right|_{\lambda = -i\tau},
\tag\labeleq{\phicur}
$$
and compute its Frenet frame
$$\aligned
 T &= \dfrac{1}{2}\Phi^{-1}\pmatrix 0 & i \\ i & 0\endpmatrix\Phi \\
 N &= \dfrac{1}{2}\Phi^{-1}\pmatrix 0 & 1 \\ -1 & 0\endpmatrix\Phi \\
 B &= \dfrac{1}{2}\Phi^{-1}\pmatrix -i & 0 \\ 0 & i\endpmatrix\Phi.\endaligned
\tag\labeleq{\fre}$$
The \BB/ formula is a consequence of the following
\proclaim{Proposition \numbit}
Let $\chi(s;\nu)$ be a solution of the linear system (\phisys )
for $\lambda = \nu$,
and let $\phi(s;\lambda)$ be a solution of (\phisys ) 
for an arbitrary $\lambda \in \C$.
Then, the gauge transformation
$$\phi^{(1)}(s;\lambda)= \pmatrix \lambda & -\nu \dfrac{\chi_1}{\chi_2} \\
-\nu\dfrac{\chi_2}{\chi_1} & \lambda \endpmatrix \phi(s;\lambda),
\tag\labeleq{\gaugesol}$$
produces a solution $\phi^{(1)}$ of
$$\dfrac{d\phi}{ds} = \dfrac{1}{2}
\pmatrix i\k^{(1)}(s) &\lambda \\ \lambda &-i\k^{(1)}(s) \endpmatrix\phi,
\tag\labeleq{\newphisys}
$$where 
$$
\k^{(1)}(s) = \k(s) + i\nu \left( \dfrac{\chi_1}{\chi_2} - 
\dfrac{\chi_2}{\chi_1} \right).
\tag\labeleq{\kone}
$$
\endproclaim
The proof is a direct verification of formulas (\gaugesol) and (\kone).  

We will
note for future reference that $z=\chi_1/ \chi_2$ satisfies the Riccati equation
$dz/ds = i \k z + \nu (1 - z^2)/2.$ 

The gauge transformation (\gaugesol) can be normalized to produce the
fundamental matrix for (\newphisys):
$$\Phi^{(1)}(s;\lambda)= \dfrac{1}{\sqrt{\lambda^2 - \nu^2}}
\pmatrix \lambda & -\nu z \\
-\nu/z & \lambda \endpmatrix \Phi(s;\lambda).
\tag\labeleq{\gauge}$$

Assume that $\k(s)$ is real, and suppose for the moment that $\lambda = -i\tau$
for $\tau \in \R$.
The transformation (\kone) does not, in general, produce a real-valued
function $\k^{(1)}(s)$ unless some condition is imposed on $\nu$.  Note that
the matrix $\Phi$ now takes values in $SU(2)$, and (\gauge) becomes
$$\Phi^{(1)}(s;-i\tau) = \dfrac{1}{\sqrt{\tau^2+\nu^2}}
\pmatrix \tau & -i\nu z \\ -i\nu/z & \tau \endpmatrix \Phi(s; -i\tau).$$
This gauge transformation is in $SU(2)$ only if
$\bar\nu \bar z = \nu/z$, which implies that $\nu$ is real and $z$ has
unit modulus.  In this case, $\k^{(1)}$ is real.
Using the reconstruction formula
$$
\tilde\gamma(s)=i {\Phi^{(1)}}^{-1}\left. \dfrac{d\Phi^{(1)}}{d\lambda}\right|_{\lambda = -i\tau},
$$
we obtain the following expression for the transformed curve:
$$\tilde{\gamma}(s) = \gamma(s) + \dfrac{i\nu}{\nu^2 + \tau^2}
\Phi^{-1}\pmatrix 0 & z \\ 1/z & 0\endpmatrix\Phi.
\tag\labeleq{\suback}$$
If we make the substitution $z=-e^{-i\beta}$, then $\beta(s)$ satisfies 
the equation (\betaODE) for the \BB/ transformation, with $\nu=C$,
and (\suback) becomes (\curveback), in terms of the Frenet frame for $\gamma$
given by (\fre).

\heading 2.2 The iterated formula \endheading
The gauge transformation (\gaugesol) can be iterated algebraically
once the solution of the linear system (\phisys) at an initial
curvature function $\k(s)$ is known. In fact, formula (\gaugesol) 
produces the eigenfunction of (\newphisys) out of which the new
gauge matrix is constructed.

We begin with some symmetry considerations. Given a solution 
$\chi(s,\nu)$ of the linear system (\phisys ) at a pair $(\k(s), \nu)$,
then 
$$
\chi(s,\bar{\nu}) = \pmatrix 0 & 1 \\ 1 & 0 \endpmatrix \overline{\chi(s,\nu)}.
\tag\labeleq{\symme}
$$
is a solution of (\phisys ) at $(\k(s),\bar\nu)$.
We use this property to describe the following iterated \BB/ transformation
$$
\gamma \buildrel \nu \over \longrightarrow
\gamma^{(1)} \buildrel \bar{\nu} \over \longrightarrow
\gamma^{(2)}.
$$
(The use of the symmetry (\symme) is crucial in constructing 
real-valued curvature functions and curves in $su(2)$.) The following
lemma is obtained by direct computation:
\proclaim{Lemma \numbit}
Let $f(s;{\bar \nu})$ be a solution of (\phisys ) at $(\k(s),{\bar \nu})$,
and let
$$
\zeta(s;\bar{\nu}) = \pmatrix \bar{\nu} & -\nu \dfrac{\chi_1}{\chi_2} \\
-\nu \dfrac{\chi_2}{\chi_1} & \bar{\nu} \endpmatrix f(s;\bar{\nu})
$$
be the new solution at $(\k^1(s), {\bar \nu})$ constructed by means of
the gauge transformation (\gaugesol ) . Then, two applications of the
gauge formula (\gaugesol), using solutions $\chi(s;\nu)$ and $\zeta(s;\bar\nu)$,
gives
$$
\split
\phi^{(2)}(s;\lambda) & = {\Cal T}(\lambda,\nu) \phi(s;\lambda)\\
& := \dfrac{1}{\sqrt{(\lambda^2 - \nu^2)(\lambda^2 - {\bar{\nu}}^2 )}}
\pmatrix \lambda^2 + |\nu|^2 \dfrac{\chi_2 \zeta_1}{\chi_1 \zeta_2}
& -\lambda ( \nu \dfrac{\chi_1}{\chi_2} + \bar{\nu} \dfrac{\zeta_1}{\zeta_2})
\\ -\lambda ( \nu \dfrac{\chi_2}{\chi_1} + \bar{\nu} \dfrac{\zeta_2}{\zeta_1})
& \lambda^2 + |\nu|^2 \dfrac{\chi_1 \zeta_2}{\chi_2 \zeta_1}
\endpmatrix \phi(s;\lambda).
\endsplit
\tag\labeleq{\due}
$$
Moreover, the matrix $\Cal{T}(-i\tau,\nu)$, $\tau \in \R$, is unitary if and
only if we choose 
$f(s,\bar{\nu})= \left( \smallmatrix 0 & 1 \\ 1 & 0 
\endsmallmatrix \right) \overline{\chi(s,\nu)}$.
\endproclaim

We can now derive the explicit formula for the transformed curve and 
for its corresponding curvature function by using the reconstruction
formula (\phicur) and the form (\due) of the new eigenfunction.
We set $\chi_1/\chi_2 = \rho e^{-i\beta}$ (no longer unit modulus),
$\lambda = -i\tau$ and we introduce the constant
$\alpha=i\dfrac{\nu^2-{\bar{\nu}}^2}
{(\tau^2 + \nu^2)(\tau^2 + {\bar{\nu}}^2)}$.
\proclaim{Proposition \numbit}
The iterated \BB/ transformation of the curve $\gamma(s)$ at the
eigenvalues $\nu$ and $\bar{\nu}$ is the constant torsion curve
$$
\gamma^{(2)}(s) = \gamma(s) + \alpha \left( pT+qN+rB \right),
\tag\labeleq{\curtwo}
$$
where
$$
\split
p & = \dfrac{-1}{(\nu-{\bar{\nu}}\rho^2)({\bar{\nu}}-\nu\rho^2)}
\left[ \dfrac{\nu -\bar{\nu}}{2i}(\tau^2+|\nu|^2)(1+\rho^2) \rho \cos\beta +
\dfrac{\nu + \bar{\nu}}{2}(\tau^2-|\nu|^2)(1-\rho^2)\rho\sin\beta \right] \\
q & = \dfrac{1}{(\nu-{\bar{\nu}}\rho^2)({\bar{\nu}}-\nu\rho^2)}
\left[ \dfrac{\nu +\bar{\nu}}{2}(\tau^2-|\nu|^2)(1-\rho^2) \rho \cos\beta -
\dfrac{\nu - \bar{\nu}}{2i}(\tau^2+|\nu|^2)(1+\rho^2)\rho\sin\beta \right] \\
r &= \tau |\nu|^2 \dfrac{\rho^4-1}
{(\nu-{\bar{\nu}}\rho^2)({\bar{\nu}}-\nu\rho^2)}.
\endsplit
$$
The associated curvature function is
$$
\k^{(2)}(s)=\k(s)+
\dfrac{2i\nu({\bar{\nu}}^2 -\nu^2) \rho}
{(\nu-{\bar{\nu}}\rho^2)({\bar{\nu}}-\nu\rho^2)}
\left[(\nu +\bar{\nu})(\rho^2-1)\cos\beta +i ({\bar{\nu}}-\nu)(\rho^2+1)
\sin\beta \right].
\tag\labeleq{\ktwo}
$$
\endproclaim
\smallpagebreak
\hedd{Remarks.}

1. We observe that the expression for $\gamma^{(2)}-\gamma$ now contains a 
binormal component proportional to the torsion of the
original curve. 

2. The double \BB/ transformation makes sense for
all $\nu \in \C$ and becomes the identity when $\nu$ is real. The formula
for the new curve depends in general on two complex parameters: $\nu$,
and the parameter $c=c_+/c_-$, assuming $\chi=c_+\chi^++c_-\chi^-$
is expressed on fixed basis $(\chi^+,\chi^-)$ of solutions of the
linear system at $(\k,\nu)$.

\heading 2.3 Closure conditions \endheading
When discussing the conditions for which the curve produced by a double
\BB/ transformation (\curtwo) is closed, we distinguish two cases:

\medskip
\hedd{Case (A)} is analogous to the situation discussed in \S 1.5 for
a single \BB/ transformation. For an initial curve $\gamma$ of length $L$,
$\chi(s;\nu)$ is taken to be an eigenfunction of
the {\it transfer matrix} $\Phi(kL,\nu)$ across $k$ periods of $\gamma$.
(However, when $\gamma$ is odd and $k$ is odd, the appropriate transfer matrix is 
$\left( \smallmatrix 0 & 1 \\ 1 & 0 \endsmallmatrix \right) \Phi(kL, \nu)$
instead.)  Then formula (\curtwo) produces a family of curves of period $kL$
parametrized by $\nu \in \C$.

However, if the initial curve is a closed elastic rod of constant torsion,
then one can show directly that, for any $\nu \in \C$, the new curve has
length $L$ and is
congruent to the original elastic rod.  To see this, first note that
the Riccati equation for $z$ is now
$$\dfrac{1}{2p}\dfrac{dz}{dx}=iz\cn x+ \dfrac{\nu}{2}(1-z^2).\tag\labeleq{\riccaz}$$
Since, when $\k=\cn x$,
$\Psi(kL;-i\sigma/2)$ has eigenvectors given by (\eigens) when $\sigma \in \R$,
similar formulas hold for the eigenvectors of $\Phi(kL;\nu)$:
$$v_+ = A
\pmatrix -\dfrac{\sigma}{2\mu}\\i\left( 1 + \dfrac{1-2\lambda_1}{2\mu}\right)
\endpmatrix,\qquad
v_- = A
\pmatrix i\left( 1 + \dfrac{1-2\lambda_1}{2\mu}\right)\\
-\dfrac{\sigma}{2\mu}\endpmatrix,\tag\labeleq{\eigens}$$
where, now,  $\sigma = 2i\nu$ for arbitrary $\nu \in \C$, and $\lambda_1$ and $\mu$ are defined in terms of
$\sigma$ as in \S 1.3.  We again observe that, since the eigenvectors are
independent of $k$, the Backlund transformation has
period $kL$ if and only if it has period $L$.  We will consider the initial
value problem for (\riccaz) corresponding to $v_+$, the other being similar;
the initial value for $z$ is then
$$z(0) = \dfrac{2\sigma - 4\mu}{p^{-2} + \sigma^2},$$
which can be rewritten as 
$$z(0)=ip\sn a - \dn a,\qquad\text{where}\ 
\dfrac{\sn a}{1+\cn a}=2p\nu.
\tag\labeleq{\zinit}$$

\proclaim{Proposition \numbit {\let\saul=\here}}
\let\saul=\here
  Let $p \in (0,1)$ be the elliptic modulus, and let $a$ lie in an open disk about
the origin in the complex plane such that $\cn(a)\ne -1$.  The solution of 
the initial value problem (\riccaz), (\zinit) is given by
$$z_+ = -(\dn x + i p \sn x )\left(\dn(x-a) + ip \sn(x-a)\right).$$
Moreover, $z_+(x)$ satisfies the identity
$$\cn x + i\nu\left(z_+ - \dfrac{1}{z_+}\right) = \cn(x-a).\tag\labeleq{\saulid}$$
\endproclaim
\demo{Proof}  For $a,x$ real, $z$ must be of unit modulus, and the substitution
$z= -e^{-i\beta}$ yields the initial value problem (\binit),(\dbdx) solved
in \S 1.5.  The solution is 
$$\beta_+(x) = -[\sin^{-1}(p\,\sn x) +\sin^{-1}(p\,\sn(x-a))]$$
and $z_+(x) = -e^{-i\beta_+}$ has the above form and satisfies the required
identity.  The rest follows by analyticity in $a$ and $x$.
\enddemo

After a single iteration of the \BB/ transformation, we have
$\k^{(1)} = \cn(x-a)$ by the above identity.  Let $y=x-a$.  
After the double \BB/ transformation, we have
$$\k^{(2)} = \k^{(1)} + i\bar\nu(w - 1/w),$$
where $w$ is a solution of 
$$\dfrac{1}{2p}\dfrac{dw}{dx}=iz\cn y+ \dfrac{\bar\nu}{2}(1-z^2).$$
The formulas in \S 2.2 show that, in performing the second \BB/ transformation,
we use the solution 
$w_+ = z_+\left(\dfrac{\bar\nu - \nu |z_+|^2}{\bar\nu|z_+|^2-\nu}\right)$.
We also know from Prop. \saul\ (with $x$ changed to $y$ and $\nu$ to $\bar\nu$)
that 
$$w=-(\dn y + i p \sn y )\left(\dn(y-\bar a) + ip \sn(y-\bar a)\right)$$
is a solution.  In fact, it can easily be verified, for example by
evaluated at $x=0$, that $w$ and $w_+$ are the same.  Now it follows by the
identity (\saulid) that $\k^{(2)} = \cn(x-a-\bar a)$, and the
new curve is congruent to the original elastic rod.
\bigskip
\hedd{Case (B).}
When the transfer matrix is a multiple of the identity
matrix, $\chi(s,\nu)$ can be taken to be an arbitrary complex linear combination
$c_+\chi^+(s,\nu) + c_-\chi^-(s,\nu)$ of solutions to (\phisys).  To see when
this is possible, let $\Delta(kL,\nu)$ be the trace of the transfer matrix;
in Floquet theory, this is known as the Floquet discriminant associated
to the linear system.  

Assume for the moment that $k$ is even or $\gamma$ is even.
Since $\det\Phi(kL,\nu)=1$, the eigenvalues
of the transfer matrix will be both 1 or both $-1$ exactly when $\Delta^2=4$.  
Moreover, if $\nu$ is a root of multiplicity two for the function $\Delta^2-4$, then the corresponding
eigenvectors will be linearly independent.  Taking into account the
 obvious changes for the case when $k$ is odd and $\gamma$ is odd,   
we summarize the discussion in the following
\proclaim{Proposition \numbit\let\floquet=\here}\let\floquet=\here 
Let $\tilde\nu$ be a complex double root of the equation
(a)  $\Delta^2(kL,\nu) - 4 = 0$, or (b) $\Delta^2(kL,\nu)+4=0$ if $\gamma$ is
odd and $k$ is odd.
Let $\chi(s)=c_+\chi^+(s;\tilde\nu) + c_-\chi^-(s;\tilde\nu)$, where $\chi^+,\chi^-$
are the columns of the fundamental matrix of (\phisys).  
Then, formula (\curtwo) produces a family of closed curves of period $nL$,
parametrized by $\omega = c_+/c_- \in\C$.
\endproclaim
We will loosely refer to the new closed curves as $k$-fold covers of the
original curve.

We should note that, when the initial curve is an elastic rod of constant torsion,
and the initial value for $\chi$ is equal to the limit,
as $\nu \to \tilde\nu$, of either $v_+$ or $v_-$, then the arguments in Case (a)
apply to show that the new curve is congruent to the original elastic rod.

In the next section we will exhibit concrete examples of new closed constant torsion 
curves obtained by double \BB/ transformations of elastic rods of constant torsion.

\heading 2.4 Exotic Curves from Elastic Rods\endheading

\topinsert
\hskip -.1in
\TrimTop{2.6in}\TrimBottom{2.6in}
\BoxedEPSF{rod13.eps scaled 300}
\hskip 0.3in
\TrimTop{2.6in}\TrimBottom{2.6in}
\BoxedEPSF{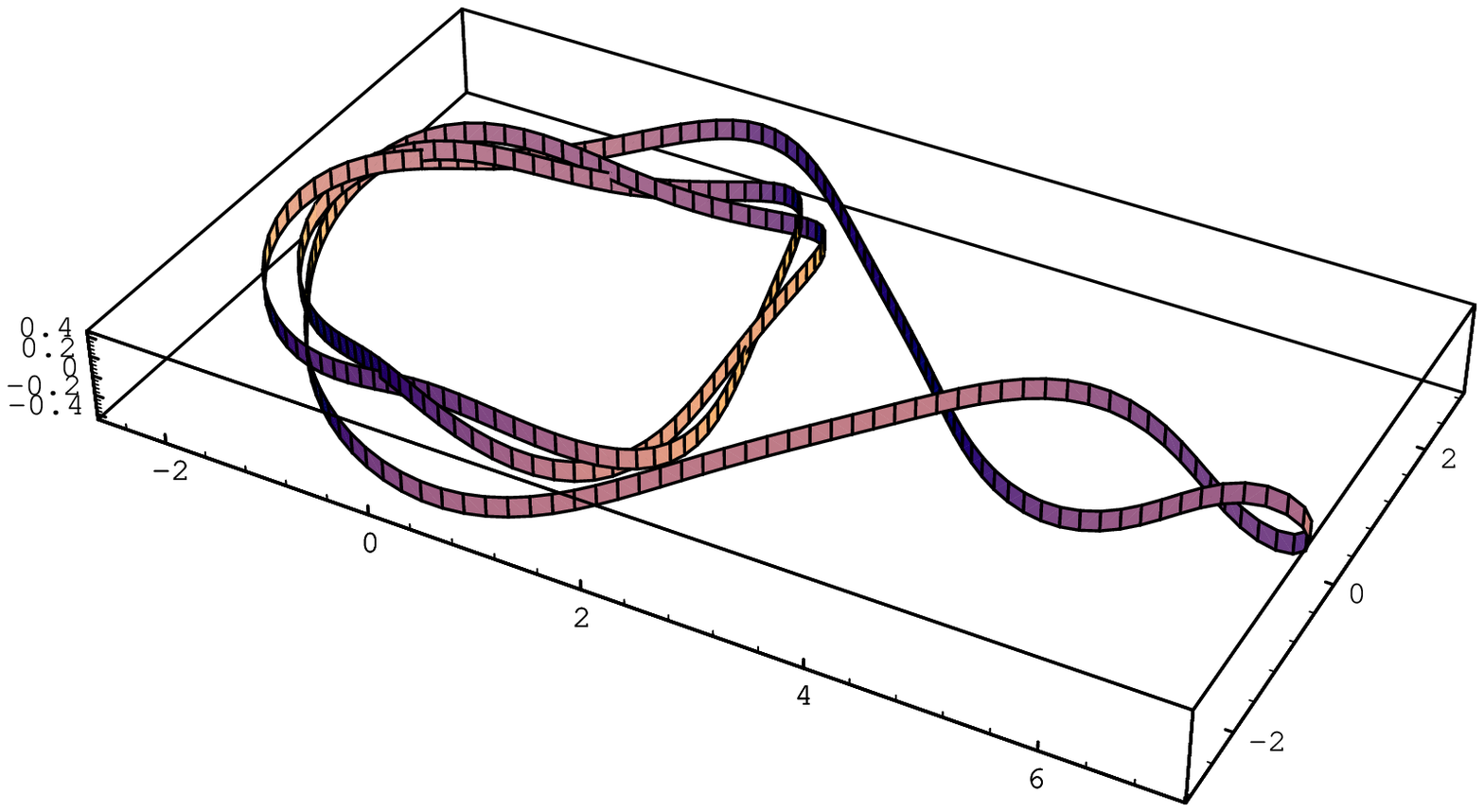 scaled 400}
\bigskip
\hskip 0.7in
\TrimTop{2.6in}\TrimBottom{2.6in}
\BoxedEPSF{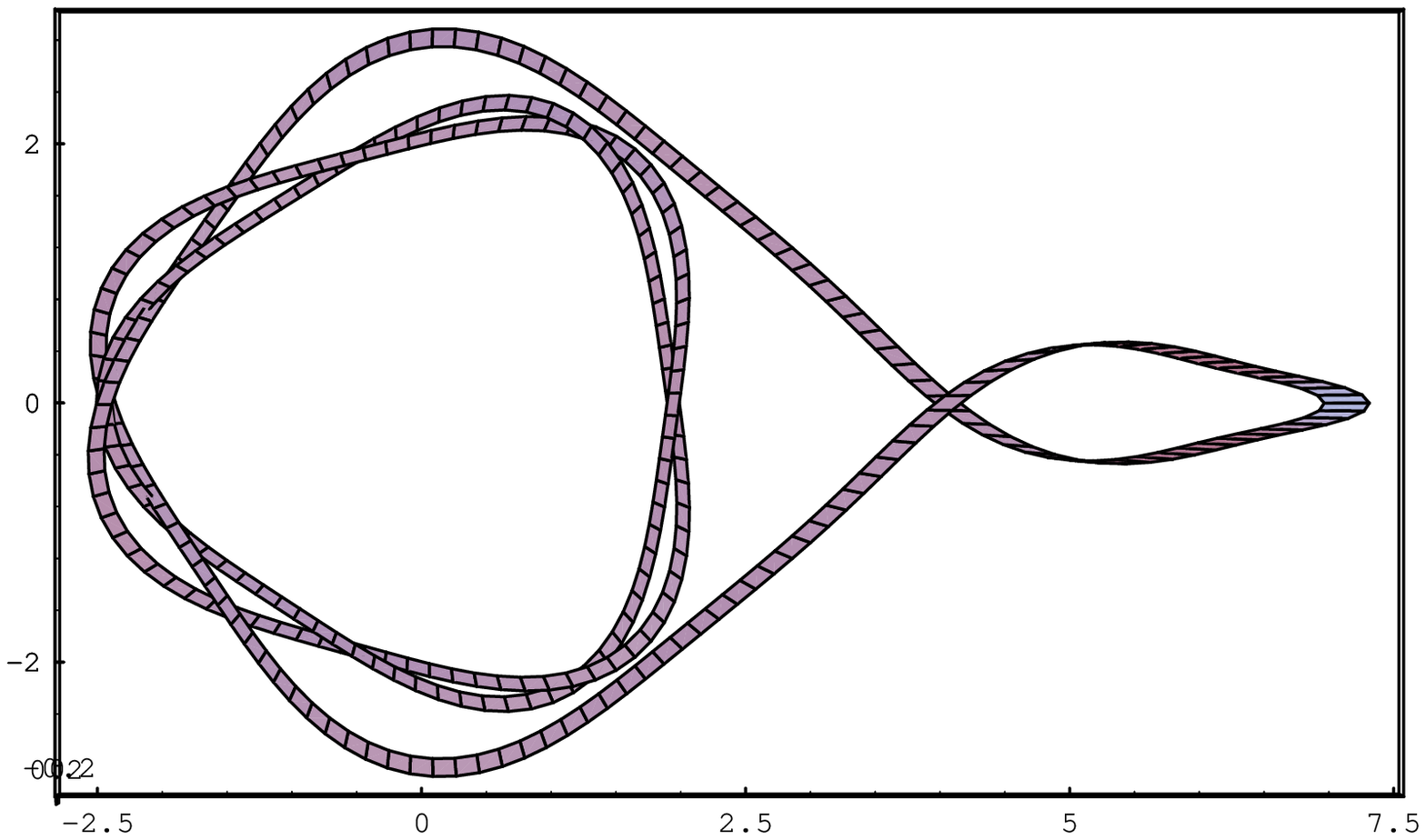 scaled 500}
\botcaption{Figure 2.1} Knotted double \BB/ transformation of the unknotted 
$(1,3)$ elastic rod
($p \approx .63093$, $n=3$, $k=4$, $\sigma \approx 1.2283 + 0.9688 i$).  In the
oblique views, the curve itself is the centerline of the vertical ribbon; the
elastic rod is at left, the transformed curve at right.  In the
view from above, the curve itself is the inner edge of the radial ribbon.
\endcaption
\endinsert

\define\tr{\text{tr}}
To produce new curves, we need to solve the equations in 
Prop. (\floquet) for
the appropriate values of $\nu$.  Suppose as before that the rod is of length $L$
and this corresponds to the elliptic parameter $x$ running from zero to $2nK$.
Then the formulas in \S 1.4  give
$$(\tr \Phi(kL))^2-4 = \tr(\Psi(kL)^2)-2=-4\sin^2(knK\Lambda)$$
when $kn$ is even, and when $kn$ is odd,
$$(\tr \left(\smallmatrix 0 &1 \\1 &0 \endsmallmatrix\right)\Phi(kL))^2+4 = 
\tr((\left(\smallmatrix 1 &0 \\0 &-1 \endsmallmatrix\right)\Psi(kL))^2)+2=
4\sin^2(knK\Lambda)$$
(Note that while the formulas in \S 1.4 apply to the linear system (\psisys), the
fundamental matrix for (\phisys) differs only by conjugation by matrix $A$, and so the trace formulas
give the same expression.)  So, in either case we look for zeros of $\sin(knK\Lambda)$, where $\Lambda$
is defined by (\lamdef) as a function of $p$ and $\sigma$, and $\sigma =2i\nu$
as above.  

In practice, we look for values of $\sigma$ in the first quadrant,
because changing $\nu$ to $\bar\nu$ does not change the double \BB/ transformation,
while changing $\nu$ to $-\nu$ is equivalent, because of symmetries in the linear
system, to just changing the sign of $\omega = c_+/c_- $.
The formula (\lambetter) shows that $\sin^2(K\Lambda)$ has simple zeros
at $\sigma = \pm p'/p \pm i$---that is, where
$\mu = \sqrt{(p^{-2} - \sigma^2)^2 + 4\sigma^2}/4$ vanishes.  (In fact these seem
to be the only complex zeros of $\sin^2(K\Lambda)$; there are, however, an
infinite number of real zeros, since $\Lambda$ is asymptotically linear in $\sigma$
when $\sigma \in \R$ (cf. formula (\lambetter)).)
We find, using Newton's method, simple zeros for $\sin(knK\Lambda)$
that have real part close to $p'/p$.  Once such a root
has been found, we can generate a family of double \BB/ transformations
parametrized by $\omega\in\C$; unless specified, the reader may
assume the value used is $\omega=1$.

We observe the following phenomena:
\medskip
\hedd{Number of roots increasing with $k$:}  
Of course, all the zeros
for $\sin(knK\Lambda)$ remain when we multiply $k$ by an integer.  However,
zeros are also picked up from $\sin(\ell K\Lambda)$ where $\ell$ is a divisor
of $kn$.  For example, when the original curve is the closed unknotted elastic
rod of constant torsion corresponding to $p\approx 0.63093$ and $n=3$
(appearing in Figure 2.1), there is
one root for $k=2$, two new roots for $k=4$, and four new roots for $k=8$.

\smallpagebreak
\hedd{Obtaining knots from unknots:}  Performing a double \BB/ transformation
of the aforementioned unknot, using one of the roots for $k=4$, yields
the knotted curve of constant torsion shown in Figure 2.1.  One can verify that
this knot has a minimum crossing number of 12, and is not a torus knot:
its Alexander polynomial is 
$$A(t)=1-t+t^3 - t^4 + t^5-t^6 + t^7 -t^9 + t^{10},$$
while the only torus knot that has an Alexander polynomial of this degree is 
a $(2,11)$-knot, whose polynomial is different.  So, this is a new knot 
type realizable by curves of constant torsion.

\smallpagebreak
\hedd{Knotting is related to parametric resonance:} 
When, for a fixed initial curve and fixed value of $k$, there are several
roots available, choosing the $\sigma$-value closest to $p'/p+i$ produces the \BB/
transform that has the most complicated shape (see below for examples).  
At $\sigma =p'/p+i$, the 
Floquet multipliers coincide but the
transfer matrix has only a single eigenvector up to multiple, and the
general solution of the linear system is unbounded. 

\smallpagebreak
\hedd{Obtaining unknots from knots:}  When the original curve is the
elastic rod of constant torsion which gives a $(2,5)$ torus knot, corresponding
to $p\approx 0.7845$ and $n=5$, three different roots are available for 
$k=2$.  Two of the resulting curves, shown in Figure 2.2, are unknotted.

\midinsert
\hskip -.1in
\TrimTop{2.6in}\TrimBottom{2.6in}
\BoxedEPSF{rod25.eps scaled 300}
\hskip 0.3in
\TrimTop{2.6in}\TrimBottom{2.6in}
\BoxedEPSF{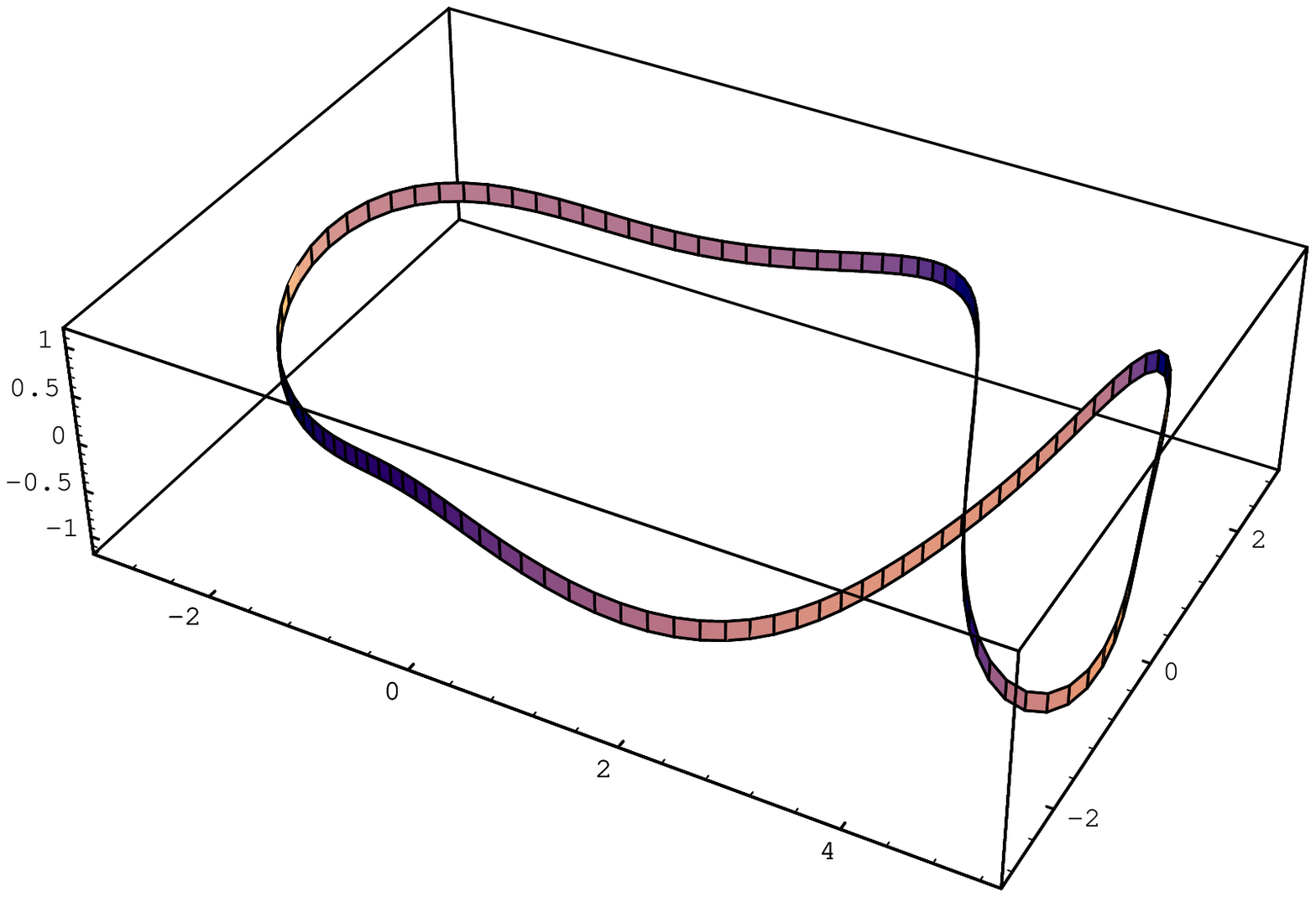 scaled 400}
\bigskip
\hskip 0.4in
\TrimTop{2.6in}\TrimBottom{2.6in}
\BoxedEPSF{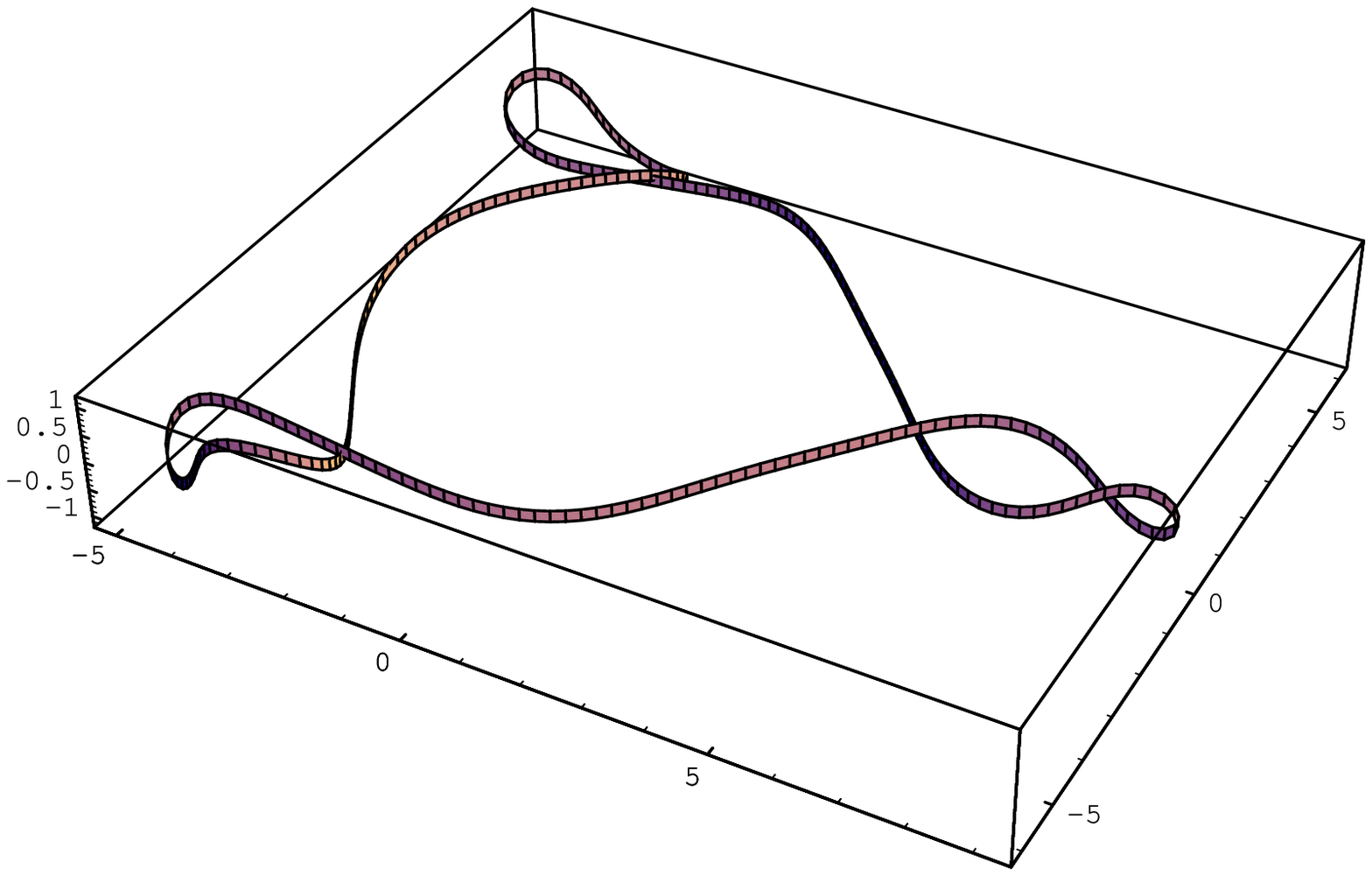 scaled 600}
\botcaption{Figure 2.2} Unknotted double \BB/ transformations of the 
knotted $(2,5)$ elastic rod, with elastic rod shown at above left
($p \approx .7845$, $n=5$, $k=2$; $\sigma \approx 0.8982 + 0.8714 i$ above right,
$\sigma \approx 0.8821 + 0.6716i$ below).
\endcaption
\endinsert
\smallpagebreak
\hedd{Self-intersections and dependence on $\omega$:}  The shape of the 
transformed curve seems to depend more on the argument of $\omega$ than its
magnitude. For $\omega$ chosen to be
real, we sometimes find that the transformed curve has self-intersections.
For example, the third available root for the $(2,5)$ knot with $k=2$ gives,
using $\omega=1$, 
the self-intersecting curve with 180-degree symmetry in Figure 2.3.  
The self-intersection, along
with the symmetry, disappears when $\omega=e^{i\pi/3}$.  The resulting
embedded curve is the twelve-crossing knot shown in Figure 2.3.  The shape of the curve changes continuously
as the argument of $\omega$ varies, showing that a given knot type can have
non-congruent realizations as a closed curve of a given constant torsion.  Moreover,
this curve has the same knot type as the
curve in Figure 2.1, showing that this knot type can be realized by closed curves of
different constant torsion! (Recall that the \BB/ transformation preserves the value
of the torsion, so these curves have the same torsion as the elastic rods we started
with.) 

\midinsert
\hskip 0.2in
\TrimTop{2.6in}\TrimBottom{2.6in}
\BoxedEPSF{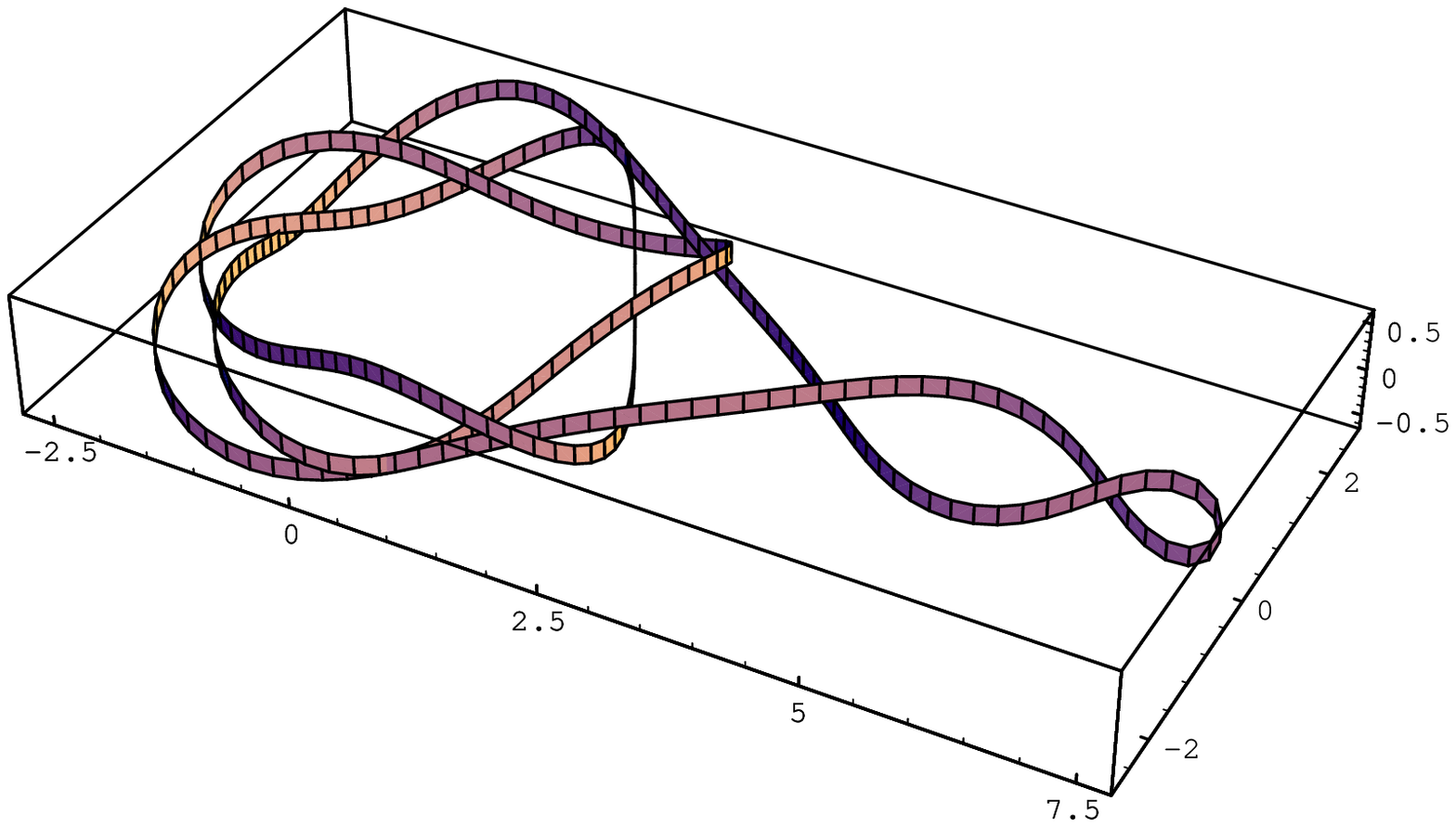 scaled 400}
\hskip 0.2in
\TrimTop{2.6in}\TrimBottom{2.6in}
\BoxedEPSF{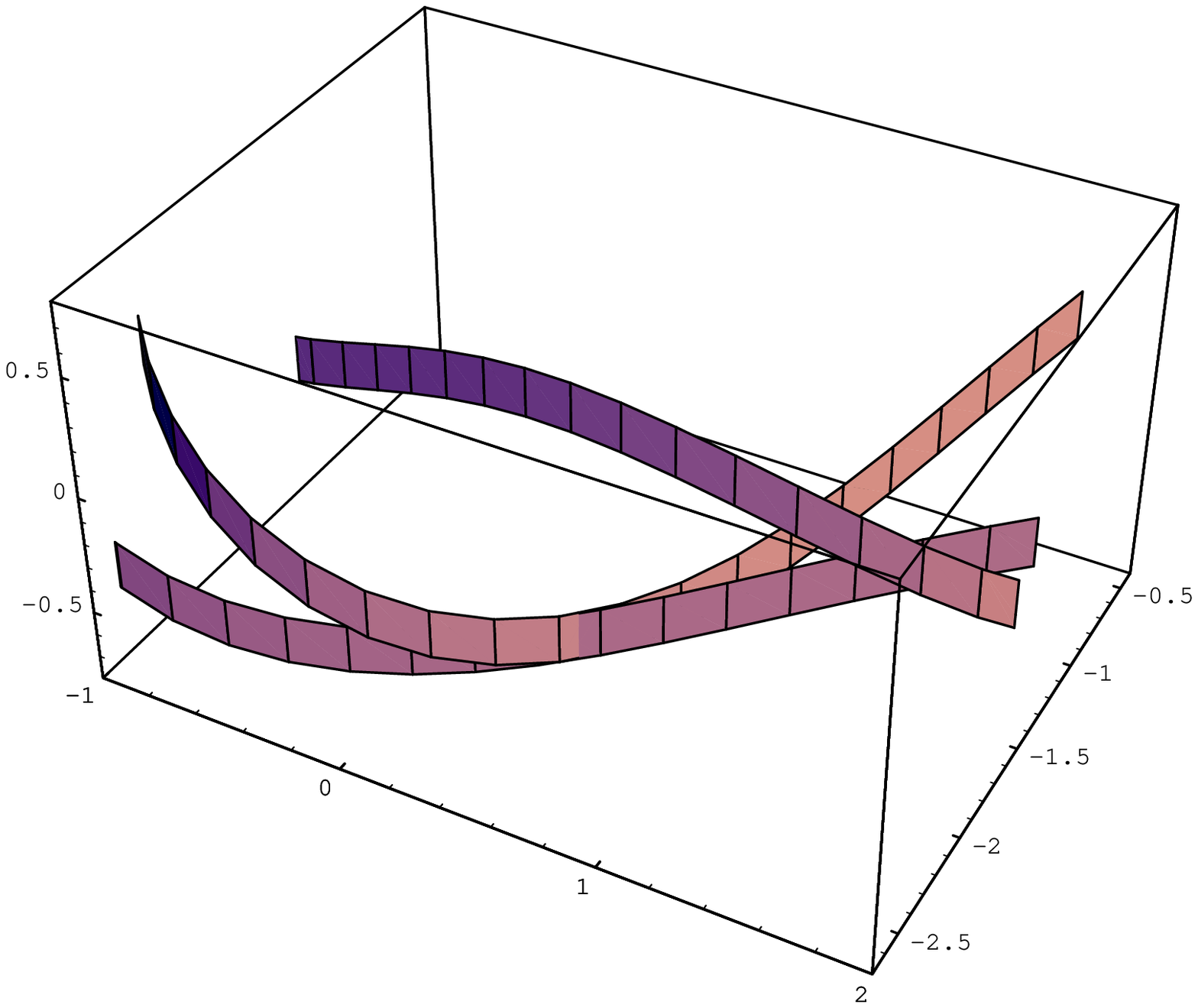 scaled 300}
\bigskip
\hskip 0.2in
\TrimTop{2in}\TrimBottom{2in}
\BoxedEPSF{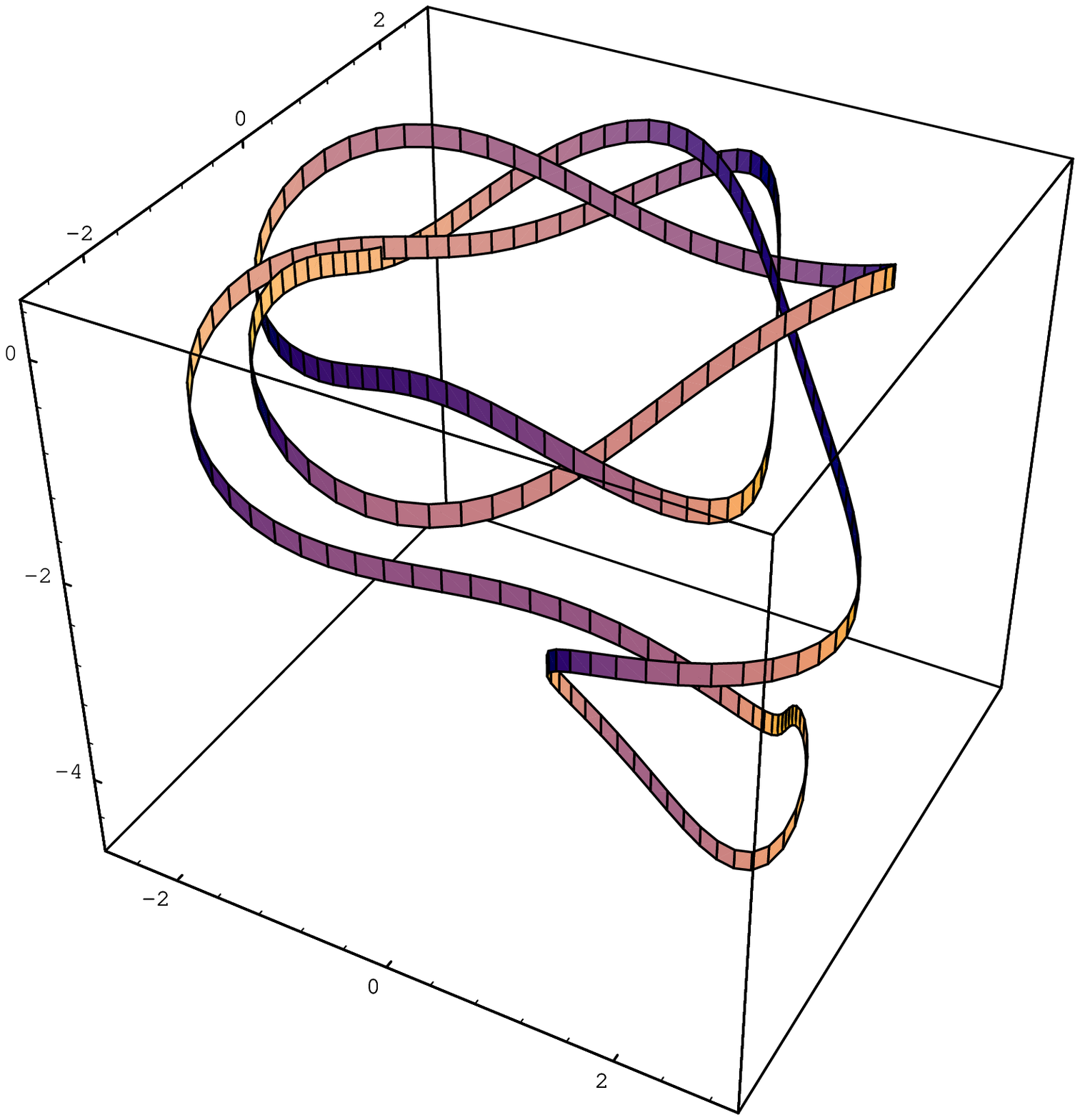 scaled 400}
\hskip 0.2in
\TrimTop{2in}\TrimBottom{2in}
\BoxedEPSF{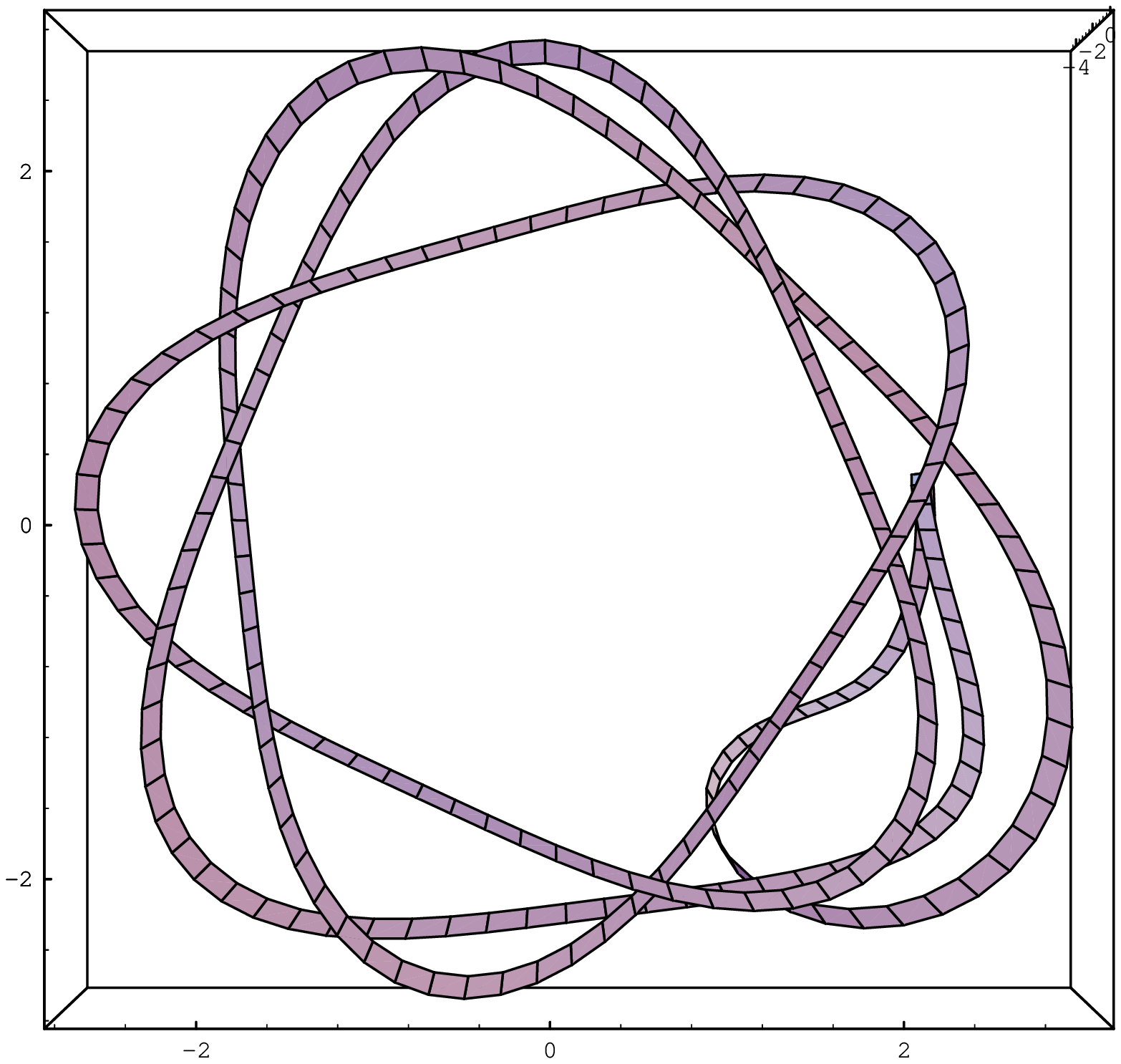 scaled 300}
\bigskip
\botcaption{Figure 2.3} Resolving self-intersection for a double \BB/ 
transformation of the $(2,5)$ elastic rod
($p \approx .7845$, $n=5$, $k=2$, $\sigma \approx 0.9067 + 0.9697 i$).  The upper
curve has two self-intersections, one of which is visible in the closeup on the
right.  The lower curve differs by a change of initial values; at right
is a view of it from above.
\endcaption
\endinsert
\bigpagebreak
Finally, we note that these exotic curves, produced by double \BB/ transformations
of closed elastic rods of constant torsion, are planar-like 3-solitons---i.e.,
a linear combination of the purely binormal LIE vector fields $B_0$,
$B_1$, $B_2$, $B_3$ vanishes along them.  This follows from our mechanical
verification, for $n=0,1,2$, of the conjecture in \S 1.7, since the
calculation is purely algebraic and applies when the \BB/ parameter $C$
takes successive complex values $\nu$ and $\bar\nu$.
\heading 2.5 Further Research \endheading
In this section, we outline some directions for future research.

\smallpagebreak
\hedd{Evolution under the LIE:}  It would be interesting to see how our exotic
curves evolve under the Localized Induction Equation:
$$\dfrac{\di\gamma}{\di t} = \k B.$$
Since $\k B$ is, up to a tangential component, a Killing field along the
elastic rods, they evolve by a rigid motion.  (In this case, the
LIE preserves constant torsion; we do not expect this to happen in general.)
In this context, Ercolani et al [E-F-M] used \BB/ transformations
of a given solution $u(x,t)$ of the sine-Gordon equations to
produce solutions that appear as homoclinic orbits, as $t$ goes from $-\infty$
to $\infty$, in the isospectral set of $u$.  For the LIE itself, the first
author [C] has produced homoclinic solutions as \BB/ transformations of
the planar circle (which translates under LIE); moreover, these solutions,
which are multiple covers of the circle, exhibit self-intersections that
persist under LIE.  It is possible that some of the self-intersecting LIE
solitons produced by our \BB/ transformations may also exhibit this behaviour.

\smallpagebreak
\hedd{The space of constant torsion knots:}  Our examples of exotic curves
of constant torsion raise the question of which knot types can be realized
with constant torsion.  For example, while the trefoil knot is not realized
by constant torsion elastic rods (it can, however, be realized by elastic
rods of non-constant torsion [I]), is it possible that by a combination of
knotting and unknotting, effected by successing double \BB/ transformations of
such elastic rods, a trefoil can be produced?  More generally,
why, in all of our examples, do all the crossings have the same sign?

\smallpagebreak
\hedd{Understanding knotting:}  The examples given in the previous section
show that, under the double \BB/ transformation, the knot type of the 
curve may change dramatically.  One approach to characterizing the change
in knot type would be to find invariants which can be computed using
the analytic representations we have for these curves.  For example, the 
M\"obius energy of the curve, which has been shown to be related to 
the minimum crossing number [F-H-W], could be computed using the Gauss integral.
It may be possible to compute the self-linking number of the new
curve and the linking number of the new curve and the old curve.  (However,
White's formula for the linking number is only valid when the ribbon stretching
between the curves is embedded; we were able to use it in \S 1.6 by taking the
\BB/ parameter $C=\nu$ to be small, but this not true for the roots $\nu$
required for the double \BB/ transformation.)  It may also be possible,
by generating and examining many more examples, to formulate conjectures 
concerning the change in knot type, in terms of constructions more familiar
to knot theorists (eg. cable knots and satellite knots).

We have shown how the Floquet spectrum of the linear system for an
elastic rod contains information about both its knot type (related
to the simple roots $\pm \frac{p}{p'} \pm i$ of the Floquet
discriminant) and the knot types of its iterated \BB/ transformations
(encoded in some way by the complex double roots). Another approach,
then, is to use Floquet analysis in combination with \BB/ transformations
for general multi-phase solutions. Exact formulas for curves of
constant torsion that generate n-phase solutions of the LIE or the
sine-Gordon equation can be constructed in terms of Theta functions 
using standard methods of algebraic geometry [Kr, D] which, at the
same time, provide a description of the associated Floquet
spectrum.  (Counting-type of arguments [P-T] can also be used to
deduce a priori information about the Floquet eigenvalues.)

\smallpagebreak
\hedd{Filling up the space of knots:}
One long-term goal is a complete classification of the knot types
of n-phase solutions. In the spirit of Gromov's knot approximation
by Legendrian curves [G], it may be possible to combine density theorems
for n-phase curves with a precise knowledge of their topology 
in order to approximate more general knotted curves.

\bigskip
\heading References \endheading
\parindent=0pt
\parskip=0pt
[B-F] P. Byrd, M. Friedman, {\it Handbook of Elliptic Integrals for Engineers
and Physicists}, Springer, 1953.

[C] A. Calini, {\it A note on a \BB/ transformation for the continuous
Heisenberg model}, Phys. Lett. A 203 (1995), 333-344.

[C-T] S.-S. Chern, C.-L. Terng, {\it An analogue of \BB/'s theorem in affine geometry},
Rocky Mountain Math J. 10 (1980), 105-124.

[D] E. Date, {\it Multi-soliton solutions and quasi-period solutions of 
nonlinear equations of sine-Gordon type},
Osaka J. Math. 19 (1982), 125-128.

[D-S] A. Doliwa, P.M. Santini, {\it
An elementary geometric characterization of the integrable motions of a curve},
Phys. Lett. A 185 (1994), 373-384.

[E] L. P. Eisenhart, {\it A Treatise on the Differential Geometry of Curves and
Surfaces}, Ginn, 1909.

[E-F-M] N. Ercolani, M. G. Forest, D. W. McLaughlin, {\it 
Geometry of the Modulational Instability, Part II: Global Results},
preprint, University of Arizona.

[F-H-W] M. Freedman, Z. He, Z. Wang, {\it M\"obius energy of knots and unknots},
Ann. Math. 139 (1994), 1-50.

[G-P] R.E. Goldstein, D.M. Petrich, {\it
Solitons, Euler's Equation, and Vortex Patch Dynamics},
Phys. Rev. Lett. 69 (1992), 555--558.

[G] M. Gromov, {\it Carnot-Caratheodory spaces seen from within},
I.H.E.S. preprint (1994).

[I] T. Ivey, {\it Knot types and homotopies of elastic rods},
in preparation.

[Kr] I.M. Krichever, {\it Methods of algebraic geometry in the theory of 
nonlinear equations}, London Math. Soc. Lecture Notes Vol. 60.
Cambridge University Press

[L-P1] J. Langer, R. Perline, {\it Poisson Geometry of the Filament Equation},
J. Nonlinear Sci. 1 (1991), 71-93.

[L-P2] ---, {\it Local geometric invariant of integrable evolution equations},
J. Math. Phys. 35 (1994), 1732-1737.

[L-S1] J. Langer, D. Singer, {\it Knotted Elastic Curves in $\R^3$},
J. London Math. Soc. (2), 30 (1984), 512-520.

[L-S2] ---, {\it Liouville integrability of geometric
variational problems}, Comm. Math. Helv. 69 (1994), 272-280.

[L-S3] ---, {\it Hamiltonian Aspects of the Kirchhoff Elastic Rod}, 
preprint (1994).

[L-S4] ---, {\it Lagrangian Aspects of the Kirchhoff Elastic Rod}, to
appear in SIAM Reviews v.38 \#4.

[M-R] H.K. Moffatt, R.L. Ricca,
{\it The helicity of a knotted vortex filament}, in 
``Topological aspects of the dynamics of fluids and plasmas'', Kluwer, 1992.

[M] D. Mumford, {\it Elastica and Computer Vision}, in 
``Algebraic geometry and its applications'', ed. C. L. Bajaj, Springer, 1994.

[Pe] R. Perline, 
{\it Localized induction and pseudospherical surfaces},
J. Phys. A 27 (1994) 5335-5344.

[Po1] W. Pohl, {\it The Self-Linking Number of a Closed Space Curve}, Journal
of Math. \& Mech. 17 (1968), 975-985.

[Po2] ---, {\it DNA and differential geometry}, Math. Intell. 3 (1980), 20-27.

[P-T] J. Poschel, E. Trubowitz, Inverse Spectral Theory. Series in 
Pure and Applied Mathematics, 130. Academic Press, 1987.

[Ri] R.L. Ricca, {\it The rediscovery of the Da Rios equations},
Nature 352 (1991), 561-562.

[Ro] C. Rogers, {\it \BB/ transformations in soliton theory}, in ``Soliton theory:
a survey of results'', ed. A. P. Fordy, St. Martin's Press, 1990.

[W] J. White, {\it Self-linking and the Gauss integral in higher dimensions},
Am. J. Math 91 (1969), 693-728.

[W-W] E.T. Whittaker, G.N. Watson, {\it A Course of Modern Analysis},
Cambridge, 1927.
\enddocument